\documentclass[12pt]{article}
\usepackage{latexsym,cite,amsthm,amsmath,amssymb}

\usepackage{amsmath}
\usepackage{amsthm}
\date{}
\def\bc{\begin{center}}
\def\ec{\end{center}}
\def\phi{\varphi}
\def\epsilon{\varepsilon}
\def\cha{{\rm char}\,}

\def\End{{\rm End}\,}
\def\Ann{{\rm Ann}\,}
\def\Reg{{\rm Reg}\,}

\def\alg{{\rm alg}}

\def\l{\lambda}
\def\g{\gamma}
\def\a{\alpha}
\def\b{\beta}
\def\d{\delta}
\def\s{\sigma}
\def\f{\phi}
\def\e{\epsilon}

\def\r{\rho}
\def\G{\it{\Gamma}}
\def\0{\bar 0}
\def\1{\bar 1}
\def\Z{\mathbb Z}
\def\ZZ{{\mathcal Z}}

\def\U{{\cal U}}
\def\L{{\cal L}}
\def\R{{\cal R}}
\def\ctd{\hfill$\Box$}
\def\HH{\mathbb H}
\def\OO{\mathbb O}
\def\QQ{\mathbb Q}
\def\bes{\begin{eqnarray*}}
\def\ees{\end{eqnarray*}}
\def\bee{\begin{eqnarray}}
\def\eee{\end{eqnarray}}

\def\le{{\langle}}
\def\re{{\rangle}}
\def\Proof{{\it Proof. }}

\newtheorem{Th}{Theorem}[section]
\newtheorem{Lem}{Lemma}[section]

\newtheorem{Cor}{Corollary}[section]
\newtheorem{Rem}{Remark}[section]

\newtheorem{Prop}{Proposition}[section]

\DeclareMathOperator{\Cay}{Cay}

\begin{document}
\title{Irreducible bimodules over alternative algebras and superalgebras}
\vspace {5 mm}

\author {Ivan Shestakov\footnote{Supported by FAPESP, Proc.\,2010/50347-9 and CNPq, Proc.\,3305344/2009-9 } {\small and} Maria Trushina\footnote{Supported by FAPESP, Proc.\,2008/50141-1}}
  \vspace {5 mm}

\maketitle
\begin{abstract}
The irreducible alternative superbimodules are studied. The complete classification is obtained for even bimodules of arbitrary dimension and for finite-dimensional  irreducible superbimodules over an algebraically closed field. 
\end{abstract}

\section{An Introduction}
\hspace{\parindent}
Alternative algebras are the nearest generalization of associative algebras. An algebra~$A$ is called alternative if it satisfies the identities
\bes
 (x, x, y)& =& 0  \ \hbox{(left alternativity),}
\\
 (x, y, y) &= &0 \ \hbox{(right alternativity),}
\ees
where $(x, y, z) = (xy)z - x(yz)$ is the associator of the elements $x$, $y$,~$z$.
 A classical example of nonassociative alternative algebra is the 8-dimensional algebra of octonions (Cayley numbers), or, more generally, a Cayley-Dickson algebra.

The structure theory of finite-dimensional alternative algebras was developed by M.\,Zorn: if $A$ is a finite-dimensional alternative algebra, then $A$ contains a unique maximal nilpotent ideal~$R$ (the nil radical of the algebra~$A$), and  the quotient algebra $A/R$ is
semisimple~\cite{Zorn2}. Every semisimple algebra is a direct sum of ideals which are simple algebras~\cite{Zorn1}.
Every simple alternative algebra is either an associative algebra
or a Cayley-Dickson algebra over its center (see\ \cite{Zorn1,Sch1}).
For principal results of the theory of infinite-dimensional alternative algebras we refer the reader to~\cite{ZSSS}.

\smallskip

The theory of representations of alternative algebras was initiated in the papers by R.D.Scha\-fer \cite{Sch2} and N.~Jacobson \cite{Jac},
 where birepresentations of finite-dimensional alternative algebras were studied.
 In particular, in~\cite{Sch2}
it has been shown that any birepresentation of a semisimple alternative algebra of characteristic~$0$
(without any restriction on dimension) is completely reducible,
and in~\cite{Jac} the irreducible alternative bimodules over finite-dimensional algebras
were classified.
It should be noted that, contrary to the associative case,
for nonassociative algebras the concept of birepresentation is more natural and easily defined than that of
 one-sided (right or left) representation.
Right representations of alternative algebras were studied
in K.~A.~Zhevlakov~\cite{Zhevl}
and A.~M.~Slin'ko, I.~P.~Shestakov~\cite{SlSh}.
In~particular, in~\cite{SlSh} some description of
right irreducible alternative modules was obtained
and a characterization of the quasi-regular radical in~terms of
representation theory was given. However, even in~the finite-dimensional case
the classification of  irreducible alternative  right modules remains
incomplete. For example, until now it is not clear whether the left regular
representation of a Cayley-Dickson algebra (it is irreducible) is an alternative right
 representation.

\smallskip

The simple alternative superalgebras were described by E.~I.~Zel'manov, I.~P.~Shestakov~\cite{ZSh}
(in~the case of characteristic different from $2$ and~$3$)
and I.~P.~Shestakov~\cite{Sh2} (for arbitrary characteristic).
It turned out that nontrivial examples of such superalgebras appear
in the case of characteristic $2$~or~$3$ only.
The study of representations of alternative superalgebras
began in N.~A.~Pisarenko~\cite{Pisar}.
In~that paper the finite-dimensional irreducible alternative
superbimodules over  finite-dimensional simple alternative superalgebras of
characteristic different from $2$~and~$3$ were described. It turned out that in~this case
nontrivial alternative non-associative superbimodules  appear
in the case of the minimal nontrivial simple superalgebra, i.e.
 doubled field $F[u] = F \oplus Fu$, $u^2 = 1$ only.
 In~\cite{ConchShest,T4} irreducible bimodules over the
simple alternative superalgebras $B(1|2)$ and $B(4|2)$
of characteristic~$3$ and of dimensions $3$~and~$6$ were described.

\medskip

The aim of our paper is to describe irreducible
birepresentations of alternative algebras of arbitrary dimension and  characteristic,
on one side, and on the other, to describe finite-dimensional irreducible alternative superbimodules.
We also describe irreducible superbimodules
of~any dimension and  characteristic over the simple alternative
superalgebras. To complete the description of irreducible
infinite-dimensional superbimodules, it remains to describe such
superbimodules (or to prove their absence) over prime locally
nilpotent alternative superalgebras of characteristic~$3$.

\smallskip

The results of the paper were announced in \cite{TSh}.

\section{Definitions and known results}

\subsection{Birepresentations of alternative algebras}

\hspace{\parindent}
Let $A$ be an alternative algebra over a field~$F$ and $V$ be
an $A$-bimodule, i.~e.\ $V$ is a vector space over~$F$
such that on $V$ left and right bilinear multiplications by the elements of the algebra~$A$
are defined:
$$
A \otimes_F V \to V,\ \ (a \otimes  v) \mapsto a \cdot v,\quad
V \otimes_F A \to V,\ \ (v \otimes  a) \mapsto v \cdot a,
$$
where $a \in A$, $v \in V$.
In the vector space direct sum
$E = A \oplus V$ define  a product~$*$ by
$$
(a_1 + v_1) * (a_2 + v_2) = a_1a_2 + a_1 \cdot v_2 + v_1 \cdot a_2,
$$
where $a_1, a_2 \in A$, $v_1,v_2 \in V$.
So $E$ becomes an algebra over~$F$,
in~which $A$ is a subalgebra and $V$ is an ideal with~zero multiplication.
The algebra~$E$ is called the \emph{ split null extension} of the algebra~$A$
by~its bimodule~$V$. The bimodule~$V$ is called an {\em alternative
$A$-bimodule}, if the split null extension $E = A\oplus V$
 is an alternative algebra. It is easy to check
that a bimodule~$V$ is alternative if and only if
for any $a, b \in A$, $v \in V$ the following equalities holds in the algebra~$E$
$$
(a, b, v) = (b, v, a) = (v, a, b),\quad (v, a, a) = 0.
$$

The description of irreducible (bi)modules is one of the main problem
in the (bi)representation theory in any class of algebras. For an algebra~$A$
denote by $\Reg A$ the \emph{regular bimodule} $V = A$
with~the action of~$A$ given by the product in~$A$.
It is clear that the bimodule $\Reg \mathbb O$ over a Cayley-Dickson algebra~$\mathbb O$  is an irreducible faithful alternative
bimodule which is not associative. Another examples of  alternative
nonassociative bimodules provide the so-called  Cayley bimodules  over a generalized quaternion algebra~$\mathbb H$, which are defined as follows.
Let $a \mapsto \bar a$ be the symplectic involution in~$\mathbb H$,
then an alternative $\HH$-bimodule $V$ is called a {\it Cayley bimodule} if it satisfies the  identity $av=v\bar a$  for any $a\in\HH,\ v\in V$. The Cayley $\HH$-bimodules form a category which is isomorphic to the category of left associative $\HH$-modules via the following isomorphism (see \cite[Lemmas 11,12]{Sh2}).
Let $L$ be a left  associative $\mathbb H$-module
with~the action $(a,v) \mapsto av$. Then the vector space $L$
 with~the bimodule operations
$$
v \cdot a = av,\quad a \cdot v = \bar av,
$$
is an alternative Cayley bimodule over $\HH$ which we will denote as $\Cay L$.


The results of the mentioned in the introduction papers of R.~D.~Schafer~\cite{Sch2}
and N.~Jacobson~\cite{Jac} lead to the following description of finite-dimensional
irreducible alternative bimodules.

\begin{Th}\label{t1}
Let $A$ be a finite-dimensional alternative algebra over a field~$F$
 of characteristic different from~$2$ and $V$ be
an irreducible alternative faithful $A$-bimodule.
Then one of the following cases holds\textup:
\begin{itemize}
\item
$A$ is associative and $V$ is an associative bimodule,
\item
$A = \mathbb O$ and $V = \Reg \mathbb O$,
\item
$A = \mathbb H$ and $V = \Cay L$, where $L$ is a minimal left ideal of $\HH$.
\end{itemize}
\end{Th}

Irreducible bimodules over  generalized quaternion algebras
were described in~\cite{Jac} under~the assumption that the ground field has
 characteristic different from~$2$ and centralizes the bimodule.
In~\cite{Sh2} the second author extended this result,
removing the restrictions on the characteristic, dimension, and the condition of centralization.

\subsection{Alternative superalgebras}\label{s2}

\hspace{\parindent}
An algebra~$A$ is called a $\mathbb Z_2$-graded
algebra or {\em superalgebra},
if $A = A_{\bar 0} \oplus A_{\bar 1}$,
where $A_iA_j \subseteq A_{i+j}$, $i, j \in \mathbb Z_2$.
The subspaces $A_{\bar 0}$ and~$A_{\bar 1}$
are called respectively the even and odd parts of the superalgebra~$A$.
For example, the Grassmann algebra $G = G_{\bar 0} \oplus G_{\bar 1}$
becomes a superalgebra, if we denote by $G_{\bar 0}$~($G_{\bar 1}$)
 the submodule generated by the words of even (odd) length
from generators of the algebra~$G$. A superalgebra
$A = A_{\bar 0} \oplus A_{\bar 1}$
is called an {\em alternative superalgebra}, if its Grassmann envelope
$G(A) = G_{\bar 0} \otimes A_{\bar 0} + G_{\bar 1} \otimes A_{\bar 1}$
is an alternative algebra.
From this definition it follows that the alternative superalgebras
are defined by the identities
\begin{align*}
& (a_i, a_j, a_k) = (-1)^{jk+1}(a_i, a_k, a_j),
\\
& (a_i, a_j, a_k) = (-1)^{ij+1}(a_j , a_i, a_k),
\\
& (a_0, a_0, A) = 0,
\end{align*}
where $a_s \in A_s$, $s = i, j, k \in \{\bar 0,\bar 1\}$.

Let $A$ be an algebra. Denote by $A[u]$
the superalgebra obtained by  doubling of the algebra~$A$:
$$
A[u] = A \oplus A \cdot u,\quad
(A[u])_{\bar 0} = A,\quad
(A[u])_{\bar 1} = A \cdot u,
$$
where $u$  is an odd central element, $u^2 = \a\in F,\ \a\neq 0$. When $u^2=1$, we will denote $A[u]$ as $A[\sqrt1\,]$.

C.~T.~C.~Wall~\cite{Wall} proved that every simple finite-dimensional
associative superalgebra over an algebraically closed field~$F$
 is isomorphic to one of the following superalgebras:
\begin{alignat*}{2}
& \bullet &\quad & A = \mathrm M_{m|n}(F),\quad
A_{\bar 0} = \left\{
\begin{pmatrix} \star & 0 \\ 0 & \star \end{pmatrix}
\begin{matrix} m \\ n \end{matrix} \right\},\quad
A_{\bar 1} = \left\{
\begin{pmatrix} 0 & \star \\ \star & 0 \end{pmatrix}
\begin{matrix} m \\ n \end{matrix} \right\},
\\
& \bullet &
& A = \mathrm M_{n}(F)[\sqrt{1}\,],
\text{  the doubled matrix algebra}.
\end{alignat*}

In~\cite{ZSh} it was proved that every nontrivial simple
alternative superalgebra of characteristic different from $2$~and~$3$
is associative. In~particular, the nontrivial
simple finite-dimensional alternative superalgebras over an algebraically closed field~$F$ of characteristic different from $2,$~$3$
are exhausted by the superalgebras
$\mathrm M_{m|n}$, $\mathrm M_n[\sqrt{1}\,]$.

\smallskip
The simple alternative superalgebras of characteristic $2$ and~$3$ were described
in~\cite{Sh2}. In~these cases nonassociative nontrivial
simple superalgebras  appeared.

Over a field of characteristic~$2$, any such a superalgebra is isomorphic to
one of the following two superalgebras obtained from a Cayley-Dickson algebra:
\begin{itemize}
\item
the Cayley-Dickson superalgebra
$\mathbb O(4|4) = \mathbb O = \mathbb H + v\mathbb H$
with the~$\mathbb Z_2$-grading
induced by the Cayley-Dickson process applying to the~generalized
quaternion subalgebra~$\mathbb H$,
\item
the doubled Cayley-Dickson algebra $\mathbb O[u]$.
\end{itemize}
Note that both these superalgebras are alternative algebras.

\smallskip

A simple nonassociative nontrivial alternative superalgebra
of characteristic~$3$ is isomorphic to one of the following superalgebras
(see \cite{Sh2, Sh3}):
\begin{itemize}
\item
\emph{superalgebra} $B(1|2)$.
Let $V$ be a two-dimensional vector space
over a field~$F$ with~a nonzero skew-symmetric bilinear form
$\langle \cdot,\cdot \rangle$.
The unital superalgebra $B(1|2) = F \cdot 1 \oplus V$
with~an~identity~$1$ is given by the grading
$$
B(1|2)_{\bar 0} = F \cdot 1,\quad
B(1|2)_{\bar 1} = V
$$
and the supercommutative multiplication
$xy = \langle x, y\rangle \cdot 1$ for $x, y \in V$;
\item
\emph{superalgebra} $B(4|2)$.
Let $V$ be the same space as above; set
$$
B(4|2)_{\bar 0} = \operatorname{End} V \cong \mathrm M_2(F),\quad
B(4|2)_{\bar 1} = V
$$
with the multiplication
$$
v \cdot a = a(v) = \bar a \cdot v,\quad
u \cdot v = \langle \cdot, u\rangle v \in \operatorname{End} V,
$$
where $a \in \operatorname{End} V$, $u, v \in V$
and $a \mapsto \bar a$ denote the symplectic
involution in~$\operatorname{End} V$
(i.~e.\ $\langle a(u), v\rangle = \langle u, \bar a(v)\rangle$
for any $u, v \in V$).
It is easy to check that $V = B(4|2)_{\bar 1}$ is a Cayley bimodule over
$\mathrm M_2(F) = B(4|2)_{\bar 0}$;
\item
\emph{twisted superalgebra of vector type} $B(\Gamma,D,\gamma)$.
Let $\Gamma$ be an associative commutative algebra over a field~$F$,
$D$ be a nonzero derivation of~$\Gamma$,
$\gamma \in \Gamma$.
Denote by $\bar\Gamma$
an isomorphic copy of the vector space~$\Gamma$ with respect to the isomorphism
$a \mapsto \bar a$ and set $ B(\Gamma,D,\gamma) = \Gamma \oplus \bar\Gamma$
with~the multiplication
\begin{align*}
& a \cdot b = ab,
\\
& a \cdot \bar b = \bar a \cdot b = \overline{ij},
\\
& \bar a \cdot \bar b = \gamma ab + 2D(a)b + aD(b),
\end{align*}
where $a, b \in \Gamma$, $ab$ is the product of the elements $a$~and~$b$
in~$\Gamma$, and with the grading
$$
B(\Gamma,D,\gamma)_{\bar 0} = \Gamma,\quad
B(\Gamma,D,\gamma)_{\bar 1} = \bar\Gamma.
$$
The superalgebra $B(\Gamma,D,\gamma)$ is simple if and only if
the algebra~$\Gamma$ does not contain proper  $D$-invariant
ideals (i.e.\ $D$-simple).
\end{itemize}

\subsection{Alternative superbimodules}\label{s3}

\hspace{\parindent}
The alternative superbimodules are defined similarly to the non-graded
case. A bimodule $V = V_{\bar 0} \oplus V_{\bar 1}$
over an  alternative superalgebra $A = A_{\bar 0} \oplus A_{\bar 1}$
is called an {\it alternative $A$-superbimodule},
if the split null extension
$A \oplus V = (A_{\bar 0} \oplus V_{\bar 0}) \oplus (A_{\bar 1} \oplus V_{\bar 1})$
is an alternative superalgebra.

We will call an $A$-superbimodule
$V^{\mathrm{op}} = V_{\bar 0}^{\mathrm{op}} + V_{\bar 1}^{\mathrm{op}}$
 {\it opposite} to an $A$-superbimodule
$V = V_{\bar 0} + V_{\bar 1}$, if
$V_{\bar 0}^{\mathrm{op}} = V_{\bar 1}$,
$V_{\bar 1}^{\mathrm{op}} = V_{\bar 0}$
and $A$ acts on it by the following rule:
$a\cdot v = (-1)^{|a|} av$, $v\cdot a = va$,
where $v \in V^{\mathrm{op}},\ a \in A_{\0}\cup A_{\1}$, and $|a|$ means the parity of $a$, that is, $|a|=i$ if $a\in A_i,\ i=\0,\1$.

It is easy to check that for any  superbimodule~$V$ the identical application $V\rightarrow V^{op},\ v\mapsto v$, is an odd isomorphism between $V$ and $V^{op}$. In particular, if $V$ is alternative, the opposite superbimodule~$V^{\mathrm{op}}$  is alternative as well. We sometimes will say that the bimodule $V^{op}$ is obtained from $V$ by {\it changing of parity}.

As in the~case of algebras, the main problem of the theory of representations of
superalgebras is the description of irreducible superbimodules. For alternative
algebras, such bimodules are always defined either over simple algebras or
over a direct sum of two simple algebras, and in the~latter case
the bimodule is associative. It is an important difference of the supercase
that there are irreducible superbimodules over
nilpotent superalgebras. This first was noted in~\cite{Sh3},
where this fact was used to construct finite-dimensional solvable
non-nilpotent alternative and Jordan superalgebras.

Let $A = A_{\bar 1} = Fx$, $x^2 = 0$.
Assume that the field~$F$ contains a primitive 3-th root
of~$1$, i.~e.\ such element~$\varepsilon$,
that $\varepsilon^2+\varepsilon+1=0$.
Consider the $A$-bimodule $V^\varepsilon(1|1) = Fv_0 \oplus Fv_1$,
where $V_{\bar 0} = Fv_0$, $V_{\bar 1} = Fv_1$
and the action of the element~$u$ is given by the equalities
$$
v_0 \cdot x = v_1,\quad
x \cdot v_0 = \varepsilon v_1,\quad
v_1 \cdot x = v_0,\quad
x \cdot v_1 = (\varepsilon + 1)v_0.
$$
It is easy to check that $V^\varepsilon(1|1)$ is
an irreducible alternative superbimodule;
furthermore, $V^\varepsilon(1|1)^{\mathrm{op}} \cong V^{\varepsilon^2}(1|1)$.

Note that the commutator bimodule
$V^\varepsilon(1|1)^{(-)}$ over the Malcev superalgebra~$A^{(-)}$
appeared in~the paper by A.~Elduque and the second author~\cite{EldSh}.

It is clear that the bimodule $V^\varepsilon(1|1)$
can be considered also as a unital
superbimodule over the superalgebra $A^{\sharp}=F \oplus Fx$ obtained from~$A$
by adjoining an identity~$1$;  obviously, it remains irreducible.

If the field~$F$ does not contain a primitive 3-th root of~$1$ and $\cha F\neq 3$,
then a four-dimensional irreducible alternative superbimodule
over~$A$ can be constructed (see~\cite{Sh3}).

\medskip
It is clear that for any simple alternative superalgebra~$A$
the regular superbimodule $\Reg A$ and its  opposite companion $(\Reg A)^{\mathrm{op}}$
are irreducible alternative superbimodules.

N.~A.~Pisarenko~\cite{Pisar} investigated the structure of alternative
superbimodules over finite-dimensional semisimple superalgebras
of characteristic different from $2$,~$3$.

We call $V$ {\it nontrivial}, if
$V_{\bar 0} \neq 0$ and $V_{\bar 1} \neq 0$.

\begin{Th}[\cite{Pisar}]\label{th2.2}
Let $A$ be a finite-dimensional semisimple alternative
superalgebra over an algebraically closed field~$F$
of characteristic different from~$2$,~$3$, thus\
$A = A_1 \oplus A_2 \oplus \ldots \oplus A_k$,
where each summand~$A_i$ is isomorphic to
either a Cayley-Dickson algebra~$\mathbb O$
or a matrix superalgebra of types $\mathrm M_{m|n}$
and $\mathrm M_k[\sqrt{1}\,]$ over~$F$. Then
\begin{itemize}
\item
any alternative superbimodule over~$A$ is a direct sum of
associative $(A_i,A_j)$-bimodules and unital alternative
$A_i$-bimodules;
\item
any nontrivial faithful irreducible alternative
superbimodule~$V$ over the simple superalgebra $A_i$ is either associative and isomorphic to one of
the bimodules $\Reg\mathrm M_{m|n}$,
$(\Reg\mathrm M_{m|n})^{\mathrm{op}}$,
$\Reg(\mathrm M_k[\sqrt{1}\,])$,
$\Reg(\mathrm M_k[\sqrt{1}\,])^{\mathrm{op}}$,
or $A_i = F[\sqrt{1}\,] = F[u]$,
$\dim V = 2$, $V_{\bar 0} = Fv_0$, $V_{\bar 1} = Fv_1$ and
the \textup(unital\textup) action of~$A_i$ on~$V$
has one of the following forms\textup:
\begin{enumerate}
\item
$v_0u = v_1$, $v_1u = 0$, $uv_0 = -2v_1$, $uv_1 = -v_0$,
\item
$v_0u = v_1$, $v_1u = 2v_0$, $uv_0 = 0$, $uv_1 = v_0$,
\item
$v_0u = \frac{1+\alpha}{\delta}v_1$,
$v_1u = \delta v_0$, $uv_0 = v_1$, $uv_1 = (1 - \alpha)v_0$,
\end{enumerate}
where $\alpha, \delta \in F$ are
non-zero roots of the equation $\alpha^2 + \alpha\delta + \delta^2 = 1$.
\end{itemize}
\end{Th}

\begin{Rem}\label{r5}
Note that in the theorem above, the bimodules of types \textup{1) - 3)}
can be defined in a unified manner as the following series of~unital bimodules over
$F[\sqrt{1}\,] = F \oplus Fu$:
$$
V = V_{\lambda,\mu}(1|1) = Fv_0 \oplus Fv_1,\quad
V_{\bar 0} = Fv_0,\quad
V_{\bar 1} = Fv_1,
$$
with the~action
$$
v_0u = (3\mu - \lambda)v_1,\quad
uv_0 = 2\lambda v_1,\quad
v_1u = 2\mu v_0,\quad
uv_1 = (\lambda + \mu)v_0,
$$
where $\lambda,\mu \in F$, $\lambda^2 + 3\mu^2 = 1$.
In fact, the bimodule of type~1) is the bimodule
$V_{-1,0}(1|1)$, the bimodule of type~2) is isomorphic to the bimodule
$V_{0,1/\sqrt{3}}(1|1)$,
the bimodules of type~3) are isomorphic to $V_{\lambda,\mu}(1|1)$ for
$\lambda = \sqrt{2-2\alpha-\delta}/2$,
$\mu = \delta/(4\lambda)$.
\end{Rem}

\medskip

The irreducible unital superbimodules over the superalgebras
$B(1|2)$ and $B(4|2)$ were described by M.~C.~L\'{o}pez-D\'{\i}az
and the second author in~\cite{ConchShest}.
(The classification of irreducible finite-dimensional alternative
 superbimodules over $B(1|2)$ was also obtained independently
 by the first author in~\cite{T4}.)

\begin{Th}[\cite{ConchShest,T4}]\label{th2.3}
Any irreducible alternative bimodule over
the superalgebra $B = B(1|2)$, which is not associative, is isomorphic either to $\Reg B$,
or to $(\Reg B)^{\mathrm{op}}$,
or belongs to the series of irreducible bimodules
\begin{gather*}
V(\lambda,\mu)(3|3) = V_{\bar 0} \oplus V_{\bar 1},
\\
V_{\bar 0} = \operatorname{Vect}_F\langle v_0,v_1R_y,v_0R_y^2\rangle,\quad
V_{\bar 1} = \operatorname{Vect}_F\langle v_1,v_0R_y,v_1R_y^2\rangle,
\end{gather*}
with the~action
\begin{align*}
& vR_y^j\cdot y =
\begin{cases}
vR_y^{j+1}, & j<2,
\\
\mu v^s, & j=2,
\end{cases}
\\
& vR_y^j\cdot x =
\lambda v^s R_y^j + jv R_y^{j-1},\quad j=0,1,2,
\end{align*}
where $v \in \{v_0, v_1\}$, $v_i^s = v_{1-i}$ and
$\lambda$,~$\mu$ are nonzero scalars.
In addition, the superbimodules $V (\lambda, \mu)$ and $V(\lambda', \mu')$
are isomorphic if and only if $(\lambda,\mu) = \pm (\lambda', \mu')$.
Moreover, $V(\lambda,\mu) \cong \bigl(V(\lambda,\mu)\bigr)^{\mathrm{op}}$.
\end{Th}

\begin{Th}[\cite{ConchShest}]\label{th2.4}
Every alternative bimodule over the superalgebra
$B(4|2)$ is completely reducible, and every irreducible bimodule up to the change
of parity is isomorphic to the regular bimodule.
\end{Th}

\subsection{Useful identities}

\hspace{\parindent}
 We recall some useful identities that hold in any alternative algebra (see \cite{BrKlein}, \cite{Sh2},\cite{ZSSS}):
\bee
[ab,c]&=&a[b,c]+[a,c]b+3(a,b,c),	\label{id2.1}\\
(ab,c,d)&=&a(b,c,d)+(a,c,d)b-(a,b,[c,d]),	\label{id2.2}\\ \
[(a,b,c),d]&=& (ab,c,d) + (bc,a,d) + (ca, b, d), \label{id2.3}\\
((c,a,b),a,b)&=&[a,b](c,a,b),\label{id2.4}\\
 ((a,b,c),x,y)&=&((a,x,y),b,c)+(a,(b,x,y),c)+(a,b,(c,x,y))\\
 &&-[b,(a,c,[x,y])]+([a,c],b,[x,y]).\label{id2.5}\\
(z,x,ty)&=&-(z,t,xy)+(z,x,y)t+(z,t,y)x,\label{id2.6}\\
	(z,x,yt)&=&-(z,t,yx)+x(z,t,y)+t(z,x,y).\label{id2.7}.
\eee
We will also need the  identities
\bee
	z(yxy)&=&((zy)x)y, \label{id2.8}\\
	(zy)(xz)&=&(z(yx))z, \label{id2.9}
\eee
known as {\it right {\rm and} central Moufang identities} (see \cite{ZSSS}).

\section{Reduction to prime superalgebras and superalgebras with $A_{\bar 1}^2=0$.}

\hspace{\parindent}
Let $A$ be an alternative superalgebra and $V$ be an irreducible $A$-super\-bi\-mo\-du\-le.
Denote by $\Ann V$ the annihilator of $V$: $\Ann V=\{a\in A\,|\, a\cdot V=V\cdot a=0\}$. It is easy to see that $\Ann V$ is an ideal of $A$ and $V$ is an irreducible faithful bimodule over $A/\Ann V.$
A bimodule $V$ is called {\it associaive} if $(V,A,A)=0$. One can easily check that if $V$ is associative than $(A,A,A)\subseteq \Ann V$, hence faithful associative bimodules exist only over associative algebras.

\smallskip

Below $V$ always will denote an irreducible faithful $A$-superbimodule which is not associative. We will use the following notations:
\begin{itemize}
\item $E=E(A,V)=A\oplus V$, the split null extension;
\item $(\lambda,\rho)$, the birepresentation of $A$ associated with $V$, that is, $\l,\r$ are linear mappings from $A$ to $\End V$ such that for any $a\in A,\ v\in V$ holds
$$
\l(a):v\mapsto (-1)^{|a||v|}a\cdot v,\ \ \ \ \r(a):v\mapsto v\cdot a;
$$
\item  $\ZZ(V)=\ZZ_A(V):=\{\phi\in\End V\,|\, [\phi,\r(a)]_s=[\phi,\l(a)]_s=0\ \forall a\in A\}$, the centralizer of the $A$-bimodule $V$, where $[a,b]_s=ab-(-1)^{|a||b|}ba$;
\item ${\mathcal E}:=\ZZ(\Reg E)$, the supercentroid of $E$;
\item $\G=\{\a\in{\mathcal E}\,|\,V\a\subseteq V,\ A\a\subseteq A\}$.
\end{itemize}
\begin{Prop}[\cite{EldSh}]\label{prop3.1}
The centralizer $\ZZ=\ZZ(V)$ is a graded division algebra and $\G$ is a commutative superalgebra which can be embedded in $\ZZ$ via the application $\pi:\a\mapsto \a|_V$.
\end{Prop}
\Proof
By the graded version of the Schur Lemma,  for an irreducible bimodule  $V$  the centralizer $\ZZ(V)$ is a graded division algebra. Consider the homomorphism
$\pi :\G\rightarrow \End V,\ \a\mapsto \a|_V$, and let $\a\in\ker \pi$. Then
$V(A\a)=(V\a)A=0$, hence $A\a=0$ and $\a=0$. It is clear also that $\pi(\G)\subseteq \ZZ$.
Furthermore, let $\a,\b\in\G,\ v\in V, x\in A$, then we have
\bes
((v\a)\b)x&=&(-1)^{|\b||x|}(v\a)(x\b)=(-1)^{|\b||x|+|\a|(|x|+|\b|)}(v(x\b))\a\\
&=&(-1)^{|\a||x|+|\a||\b|}((v\b)x)\a=(-1)^{|\a||\b|}((v\b)\a)x.
\ees
Therefore, $(V[\a,\b]_s)A=0$, which implies that $[\a,\b]_s=0$ and $\G$ is a supercommutative subsuperalgebra of $\ZZ$.

\ctd

\begin{Th}\label{th3.1}
Let $A=A_{\0}\oplus A_{\1}$ be an alternative superalgebra. If there exists an irreducible faithful $A$-superbimodule
 $V$ which is not associative  then $A$ is prime or $A_{\1}\neq 0,\, A_{\1}^2=0.$
\end{Th}

We first prove two lemmas.

\begin{Lem} \label{lem3.1}
Let nonzero elements $a,b\in A_{\0}\cup A_{\1}$ satisfy $[a,b]_s=(a,A,b)=0$. Then the subset $(V,a,b)$ is an $A$-subbimodule of $V$. Moreover, if $ab=ba=0$ and at least one of the elements $a,b$ is even  than $(V,a,b)=0$.
\end{Lem}
\proof  In fact, for any $v\in V_{\0}\cup V_{\1},\, r\in A_{ \0}\cup A_{\1}$ we have  by (\ref{id2.2}) and its superization
\begin{eqnarray*}
(v,a,b)r&=&\pm(vr,a,b)\pm v(r,a,b)\pm (v,r,[a,b]_s)=\pm (vr,a,b)\in (V,a,b),\\
r(v,a,b)&=&(rv,a,b)\pm (r,a,b)v+ (r,v,[a,b]_s)=(rv,a,b)\in (V,a,b),
\end{eqnarray*}
which proves that $(V,a,b)$ is a subbimodule of $V$.  Now, let $ab=ba=0$ and $a\in A_{\0}$. If $(V,a,b)\neq 0$ than by irreducibility $(V,a,b)=V$. By the Moufang identities, we have $a(v,a,b)=(v,a,ba)=0,\ (v,a,b)a=(v,a,ab)=0$. Therefore, $0\neq a\in  \Ann V$, a contradiction.
\ctd

\begin{Lem} \label{lem3.2}
Let nonzero ideals $I,J$ of $A$ satisfy $IJ=JI=0$. Then  $I+J\subseteq A_{\1}$. Moreover, if $A=A_{\0}$ then $A$ is prime.
\end{Lem}
\proof Assume that $I_{\0}\neq 0$. By the previous lemma,  the set $N=(V,I,J)$ is  an $A$-subbimodule of $V$ and $(V,I_{\0},J)=(V,I,J_{\0})=0$. If $i\in I_{\1},\, j\in J_{\1},\, a\in I_{\0},$ then, by the Moufang identities, for any homogeneous $v\in V$
\bes
 (v,i,j)a &=& (v,i,aj)+(v,a,ij)+(v,a,j)i = 0,\\
a(v,i,j) &=& (v,i,ja)-(v,a,ji)-(-1)^{|v|}i(v,a,j) = 0 .
\ees
Thus $0\neq I_{\0}\subseteq \Ann N$, and $N=0$. Consider now the sets $I\cdot V,\ V\cdot I$. It is easy to see that they are subbimodules of $V$ and at least one of them is nonzero. Assume that $I\cdot V\neq 0$, then $I\cdot V=V$ and we have
$J\cdot V=J\cdot(I\cdot V)\subseteq (JI)\cdot V+(J,I,V)=0$. Similarly, if $V\cdot I\neq 0$ then we would have $V\cdot J=0$ and $0\neq J\subseteq \Ann V=0$, a contradiction. Therefore, $V\cdot I=0$, and consequently $(V,A,I)=0$. Furthermore,  by the Moufang identity we have
$$
I\cdot (V,A,A)\subseteq (V,A,AI)+(V,I,AA)+A\cdot(V,I,A)=0,
$$
which proves that $(V,A,A)\subseteq \Ann_V(I)=\{v\in V\,|\,v\cdot I=I\cdot v=0\}$.
It is clear that $\Ann_V(I)$ is a subbimodule of $V$ which should be zero since $I\neq 0$. Therefore, $(V,A,A)=0$ which contradicts to non-associativity of $V$.
The contradiction proves that $I_{\0}=0$. Similarly, $J_{\0}=0$. In particular, if $I$ is an ideal of $A$ with $I^2=0$ then $I\subseteq A_{\1}$.

\smallskip

We continue by considering the case $A=A_{\0}$. Assume that $I,J$ are ideals of $A$ with $IJ=0$. Then $(I\cap J)^2=0,\ I\cap J\subseteq A_{\1}=0$, and we have $JI\subseteq I\cap J=0$. Therefore, $I=0$ or $J=0$, and $A$ is prime.
\ctd

\smallskip

 {\it The proof of the theorem.} Assume that $A_{\1}^2\neq 0$ and prove that in this case $A$ is prime. As above, it suffices to prove that if $I,J$ are ideals in $A$ with $IJ=JI=0$ then $I=0$ or $J=0$. Assume that   $I\neq 0,\, J\neq 0$, then by Lemma~\ref{lem3.2}, we have $I+J\subseteq A_{\1}$. Therefore, $IA_{\1}+A_{\1}I\subseteq I_{\0}=0$ and $0\neq A_{\1}\subseteq \Ann I$.
 Since the annihilator of an ideal in $A$ is an ideal, we may apply Lemma \ref{lem3.2} to the ideals $I, \Ann I$ and to get $\Ann I\subseteq A_{\1}$. Then $A_{\1}^2\subseteq \Ann I\subseteq A_{\1}$ and $A_{\1}^2=0$, a contradiction.
\ctd


\section{The case $A_{\1}=0.$} \label{A=A_0}

\hspace{\parindent}
The objective of this section is to generalize the results by N.Jacobson and R.D.Schafer on irreducible alternative bimodules given in Theorem \ref{t1} to arbitrary dimension and characteristic.
\begin{Th}\label{th4.1}
Let $A$ be an alternative algebra and $V$ be an irreducible
alternative faithful   $A$-bimodule which is not associative, both $A$ and $V$ to be of arbitrary dimension. Then one of the following cases holds\textup:
\begin{itemize}
\item
$A = \mathbb O$ and $V = \Reg \mathbb O$,
\item
$A = \mathbb H$ and $V = \Cay L$, where $L$ is a minimal left ideal of $\HH$.
\end{itemize}
\end{Th}

\proof
By Theorem \ref{th3.1}, $A$ is prime. Hence, by \cite{ZSSS}, $A$ is either an associative prime algebra or  a Cayley-Dickson ring or is degenerate, that is, contains nonzero absolute zero divisors.
(Recall that an element $a\in A$ is called an {\it absolute zero divisor (a.z.d.)} if $aAa=0$).

In the last case, let $0\neq a\in A$ be an a.z.d. We next show that $a\in \Ann V.$ Indeed, let $u=v\cdot a\neq 0$ for some $v\in V$. Then $V$ is generated by $u$.
Let $(\lambda,\rho):A\rightarrow \End V$ be the birepresentation of $A$ associated with the $A$-bimodule $V$ and  $M(A)$ be the subalgebra of the associative algebra $\End V$ generated by the set $\{ \lambda(a),\rho(a)\,|\,a\in A\}$.
Then $V=uM(A)$, and  there exists $W\in M(A)$  such that $v=uW=v\r(a) W=\cdots = v(\r(a) W)^n$.
 From \cite[Corollary of Lemma 1, Corollary 1 of Theorem 1]{Sh1},  it follows that the operator $\rho(a)W$ is nilpotent. Hence $v=0$, a contradiction.

\smallskip

Therefore, $A$ is a prime associative algebra or  a Cayley-Dickson ring. Let us prove that in the first case $A$ is a central order in a generalized quaternion algebra.
Show first that $A$ is not commutative. In fact, if it were true then by Lemma \ref{lem3.1} for any $a,b\in A$ we would have an $A$-subbimodule $(V,a,b)$. Since $(V,A,A)\neq 0$, there exist $a,b\in A$ such that $V=(V,a,b)$. But then $V=(V,a,b)=((V,a,b),a,b)=[a,b](V,a,b)=0$ by (\ref{id2.4}), a contradiction.

\smallskip

Let $V_{as}=\{v\in V\,|\, (v,A,A)=0\}$. Since $A$ is associative, it follows from (\ref{id2.2}) that $V_{as}$ is a subbimodule of $V$ which should be zero since $V$ is irreducible and not associative.
Since $A$ is prime and noncommutative, $A$ does not satisfy the identity $[x,y]^4=0$. Let $a,b\in A$ such that $n=[a,b]^4\neq 0$. Assume that there exists $c\in A$ such that $[n,c]\neq 0$. By \cite[Lemma 7.5]{ZSSS}, we have $((V[n,c],A,A),A,A)=0$.
Therefore, $((V[n,c],A,A)\subseteq V_{as}=0$ and $V[n,c]\subseteq V_{as}=0$. Similarly, $[n,c]V=0$ and $[n,c]\subseteq \Ann V=0$. The contradiction shows that $A$ satisfies the identity $[[x,y]^4,z]=0$. Since $A$ is prime,  the Posner theorem \cite[x]{Rowen} implies that $A$ is a central order in a generalized quaternion algebra.

\smallskip

Therefore,  $Z(A)\neq 0$ and the algebra of quotients $(Z^*)^{-1}A$ is a generalized quaternion algebra $\HH$ or a Cayley-Dickson algebra $\OO$ over $(Z^*)^{-1}Z.$
Let us show that $V$ is a $(Z^*)^{-1}A-$bimodule.

\begin{Lem}\label{lem4.1}
Let $A$ be a prime superalgebra, $Z=Z(A)_{\0}\neq 0$. Then $\l(z)=\r(z)\in\ZZ(V)$ for any $z\in Z$ and $V$ is an irreducible faithful $(Z^*)^{-1}A$-bimodule.
\end{Lem}
\Proof
By Lemma \ref{lem3.1}, for any $z\in Z,\ a\in A$ the set $(V,a,z)$ is an $A$-subbimodule of $V$. If $(V,a,z)=V,$ then, by (\ref{id2.4}) we have $V=(V,a,z)=((V,a,z),a,z)=0,$ a contradiction.  Hence $(V,a,z)=0$, and we have
  	\[(Vz)a=V(za)=V(az)=(Va)z\subseteq Vz,\]
 		\[a(Vz)=(aV)z\subseteq Vz,\]
which proves that	$Vz$ is an $A$-subbimodule of $V.$ Since $V$ is faithful, we have $Vz=V$  for any nonzero $z\in Z$.  		Moreover, $\Ann z$ is an $A$-subbimodule of $V$ too. Hence $\Ann z=0.$ In other words, the aplication $\r(z)\in\End V$ is inversible.

 		We proceed to show that $[V,z]=0.$ By (\ref{id2.1}), the set $[V,z]$ is an $A$-subbimodule of $V.$ If $[V,z]=V,$ then, by (\ref{id2.2}), for any $a,b\in A$
 		$$
		(V,a,b)=([V,z],a,b)\subseteq (ab,V,z)+a(b,V,z)+(a,V,z)b=0,
		$$
 a contradiction.   Thus $[V,z]=0,\ \r(z)=\l(z)\in\ZZ(V)$ and $V$ has a structure of a vector space over the field $(Z^*)^{-1}Z.$ Clearly, $V$ has also a natural structure of a  $(Z^*)^{-1}A-$bimodule, which is evidently  faithful and irreducible.

	\ctd
\smallskip

   If $(Z^*)^{-1}A=\OO$,  a Cayley-Dickson algebra, then $V\cong \Reg\OO$,  the regular bimodule  \cite{Jac}.   Consider an element $m\in V$ corresponding to the identity 1 of $\OO$; then $V=m\cdot\OO.$ On the other hand,  $m\cdot A$ is a nonzero $A$-subbimodule of $V.$ Hence $m\cdot A=V=m\cdot\OO.$ Consequently, for every $z\in Z$ there exists $a\in A$ such that $m\cdot z^{-1}=m\cdot a.$ We have $m\cdot (1-az)=0.$ Therefore $1=az$ and $z^{-1}=a.$ Thus $A=\OO$ is a Cayley-Dickson algebra.

\smallskip

 If $(Z^*)^{-1}A=\HH$,  a generalized quaternion algebra, then $V$ is the Cayley bimodule $\Cay L$, that is, $V$ is isomorphic to the left irreducible associative $\HH$-module $L$ on which $\HH$ acts as follows:  $v\cdot a=av$ and $a\cdot v=\bar{a}v$, where $av$ is the associative action of $\HH$ on $L$ and $a\mapsto \bar{a}$ is the symplectic involution in $\HH$.  In this case, 		consider $I=\Ann_{l} \,L=\{a\in A\,|\, aL=0\}$. Obviously, $I$ is an ideal of $A$ and $V\cdot I=0.$ Hence $(A,V,I)=0.$ Using the Moufang identities, we obtain that $(V,A,A)I=I(V,A,A)=0.$ Note that, by (\ref{id2.2}), $(V,A,A)$ is an $A-$subbimodule of $V.$ It follows that $I=0,$ because $V$ is not associative. Thus $L$ is a faithful irreducible left $A$-module. Therefore, $A$ is a primitive $PI$-algebra and, since $A$ can be embedded into $\HH$, $A$ is a $PI$-algebra. By Kaplansky's Theorem, it follows that the center of $A$ is a field. Hence $A=(Z^*)^{-1}A=\HH$.

\ctd
 		
\section{The case $A_{\1}^2=0, A_{\1}\neq 0.$}
\subsection{$A=F\,1\oplus F\, x,\ x^2=\l.$}
\hspace{\parindent}
We will first consider the minimal possible case when $\dim_FA_{\1}=1,\ A_{\1}=Fx,\,x^2=0$. The category of alternative bimodules over the superalgebra $A=A_{\1}=Fx$ is isomorphic to the category of unital alternative bimodules over the superalgebra $A^{\sharp}=F1\oplus Fx$, therefore we will study the last one.
In fact, in order to include also Pisarenko's result on bimodules over $F[\sqrt{1}]$, we will consider the more general case when
$A=F\,1\oplus F\, x,\ x^2=\l\in F.$

Recall the definition of the universal multiplicative enveloping superalgebra $\U(A)$  of an alternative superalgebra $A$.  Let $\bar A$ be an isomorphic copy of the vector space $A$ under the isomorphism $a\mapsto \bar a$. Consider the  vector space direct sum $A\oplus\bar A$ with the $\Z_2$-grading induced by that on $A$. Then the
tensor algebra $T(A\oplus \bar A)$ has a natural structure of an associative superalgebra. Denote by $I_{alt}$ the ideal of  $T(A\oplus \bar A)$ generated by the set of elements
\bee
&& a\otimes b-ab+(-1)^{|a||b|}(b\otimes a-ba),\nonumber \\
&&\bar a\otimes \bar b-\overline{ab}+(-1)^{|a||b|}(\bar b\otimes \bar a-\overline{ba}),\nonumber \\
&& ab-a\otimes b-\bar a\otimes b+(-1)^{|a||b|} b\otimes \bar a,   \label{Ialt}\\
 &&\overline{ab}+\bar a\otimes b-(-1)^{|a||b|}(\bar b\otimes\bar a+ b\otimes \bar a),\nonumber \\
&&c\otimes c-c^2,\nonumber
 \eee
where $a,b\in A_{\0}\cup A_{\1},\ c\in A_{\0}$; the  elements of the last type are needed only in case of $\cha F=2$. Since the generators of the ideal $I_{alt}$ are homogeneous, the quotient algebra $\mathcal{U}(A)=T(A\oplus\bar A)/I_{alt}$ inherits the superalgebra structure from $T(A)$.

Consider a pair of linear mapping $(\L,\R):A\rightarrow \U(A)$ defined as follows:
$\L:a\mapsto \L_a:=\bar a+I_{alt},\ \R:a\mapsto \R_a:= a+I_{alt}$; then for any alternative birepresentation $(\l,\r):A\rightarrow \End V$ there exists a unique  (super)algebra homomorphism $\phi:\U(A)\rightarrow \End V$ such that $\l=\phi\circ\L,\ \r=\phi\circ\R$. Conversely, every homomorphism $\phi:\U(A)\rightarrow\End V$ defines a structure of an alternative $A$-bimodule on $V$. In other words, the categories ${\rm Bimod_{Alt}}$--$A$ and ${\rm Mod_{As}}$--$\U(A)$ are isomorphic. The (super)algebra $\U(A)$ is called {\it the universal multiplicative  enveloping (super)algebra} of $A$.

If $A$ has a unit 1 then the quotient algebra $\U_1(A)=\U(A)/(1-\R_1,1-\L_1)$ is called {\it the universal unital  multiplicative  enveloping (super)algebra} of $A$.

\begin{Prop}\label{prop5.1} Let $A=F1+Fx,\ x^2=\l,\ A_{\0}=F1,\ A_{\1}=Fx$. Then the universal unital multiplicative enveloping superalgebra $\U_1(A)$ is a free $4$-dimensional module over its center which is isomorphic to the polynomial ring $F[t]$. More exactly, it has a basis   $\{1,a,b,ab\}$ over  $F[t]$, where $a^2=t,\ b^2=t-2\l,\ a\circ b=\l-t$.

If $\cha F=3$ and $\l=0$ then $\U_1(A)$ has a graded ideal $I=F[t](ab-t-\l)+F[t](a-b)$ such that $I^2=0$, and in this case $\U_1(A)/I\cong (F[t])[\sqrt t]$.

If $\cha F\neq 3$ or $\l\neq 0$ then $\U_1(A)$ is a prime superalgebra whose central closure is a generalized quaternion algebra over the field $F(t)$.
  \end{Prop}
\Proof Denote $a=\R_x,\ b=\L_x$, then $\U_1(A)=\alg\le 1,a,b\re$. It follows from (\ref{Ialt}) that the elements $a,b$ satisfy the relations
\bes
a^2-\l=b^2+\l=-a\circ b.
\ees
Denote $t=a^2$, then $t$ is an even central element of $\U_1(A)$, and one can easily check that $\U_1(A)$ is spanned as an $F[t]$-module by the elements $1,a,b,ab$. Moreover, a standart use of the Groebner-Shirshov basis method \cite{BK} (known also as The Diamond Lemma) shows  that  $\U_1(A)$ is a free $F[t]$-module and $F[t]$ is a polynomial ring.

If $\cha F=3$ and $\l=0$ then one can easily check that $I$ is an ideal in $\U_1(A)$ and $I^2=0$. The quotient superalgebra $\bar\U=\U_1(A)/I$ has the form $\bar\U=F[t]\oplus F[t]\,\bar a,\ \bar a^2=t$, where  $\bar\U_{\0}=F[t],\ \bar\U_{\1}=F[t]\,\bar a$.

Now let $\cha F\neq 3$ or $\l\neq 0$. Since $\U_1(A)$ is a free module over the center $F[t]$, we may consider the central closure $\tilde\U=(F[t])^{-1}\U_1(A)$ which has the same basis $\{1,a,b,ab\}$ over the field $F(t)$. Denote $u=-\tfrac{1}{t-\l} ab$, then $u^2=u-\tfrac{t(t-2\l)}{(t-\l)^2}$ and the subalgebra $F(t)[u]$ has a structure of a composition algebra over $F(t)$ with the involution $u\mapsto \bar u=1-u$. One can easily check that $\tilde\U\cong (F(t)[u],t)$, that is, $\tilde\U$ is obtained from $F(t)[u]$ by the Cayley-Dickson process with the parameter $t$ \cite[section 2.2]{ZSSS}. Hence $\tilde\U$ is a generalized quaternion algebra.

\ctd

We can now classify irreducible bimodules over the superalgebra $A=F1+Fx, \ x^2=\l$.

\begin{Th}\label{th5.1}
Let $V$ be a unital irreducible faithful alternative non-associative bimodule over the superalgebra $A=F1+Fx,\ x^2=\l$. Then there exists a simple algebraic field  extension $K=F(\a),\ \a\neq \l$, such that $V$ is a $K$-vector space of $\dim_KV\leq 4$, and up to the changing of parity  we have the following possibilities.

If  the polynomial $f(u)=u^2+(\a-\l) u+\a(\a-2\l)$ has  a root  $\epsilon\neq 0$ in $K$,    then $V=Kv\oplus Kvx,\ V_{\0}=Kv,\ V_{\1}=Kvx$, and the action of $A$ is given as follows:
\begin{itemize}
\item[$V_{\a,\l}^{\epsilon}(1|1):$] if $\a\neq 0$ then
\bes
&v\cdot x=vx,& x\cdot v=\tfrac{\epsilon}{\a}\, vx,\\
&vx\cdot x=\a v,&  x\cdot vx=(\a+\epsilon-\l)\,v;
\ees
\item[$V_{\l}(1|1):$]  if $\a=0$ then
\bes
&v\cdot x=vx,& x\cdot v=-2 vx,\\
&vx\cdot x=0,&  x\cdot vx=-\l\,v.
\ees
\end{itemize}

If the polynomial $f(u)$ above is irreducible over $K$, then $\U(\a):=\U_1(A)/(t-\a)$ is a graded division algebra over $K$, and $V$ is isomorphic to the left regular super-module  over $\U(\a)$. More exactly,
\begin{itemize}
\item [$V_{\a,\l}(2|2):$]
$\dim_KV=4,\ V$ has a basis $\{v,\,xv,\,vx,\, x(vx)\}$, $V_{\0}=Kv\oplus Kx(vx)$, $V_{\1}=Kvx\oplus Kxv$, and the action of $A$ is given as follows:
\bes
v\cdot x=vx,&& x\cdot v=xv,\\
vx\cdot x=\a v,&&  x\cdot vx=x(vx),\\
xv\cdot x=x(vx)+(\l-\a)v,&& x\cdot xv=(2\l-\a) v,\\
 x(vx)\cdot x=(\a-\l)vx+\a xv,&& x\cdot x(vx)=(2\l-\a)vx.
   \ees
 \end{itemize}
Conversely, for any simple algebraic field extension $K=F(\a),\ \a\neq\l$, the above defined  superbimodules $V_{\a,\l}^{\epsilon}(1|1)$ and $V_{\a,\l}(2|2)$ are irreducible and alternative.
\end{Th}
\Proof
Notice first that every irreducible alternative bimodule $V$ over a unital (super)algebra is necessary unital. In fact, by Lemma \ref{lem4.1}, $\l(1)=\r(1)\in \ZZ(V)$, which implies that $W=\{v-v\cdot 1\,|\,v\in V\}$ is a subbimodule of $V$. Moreover, $W\neq V$ since $1\in\Ann W$ and $\Ann V=0$. Thus $W=0$, and $V$ is unital.

Therefore, the bimodule $V$ in the theorem is a unital associative irreducible $\U$-module for $\U=\U_1(A)$. Let $\phi:\U\rightarrow \End V$ be the representation of $\U$ corresponding to the $\U$-supermodule $V$, then $\phi(t)=\a$ lies in the even part of the centralizer of $V$ which is a graded division algebra. Thus we may consider $\phi(\U)\subseteq \End V$ as a homomorphic image of the superalgebra $\U(\a)=\U/(t-\a)$ over the field $K=F(\a)$ with the same basis $\{1,a,b,ab\}$, where $a^2=\a,\ b^2=\a-2\l,\  a\circ b=\l-\a$. It is easy to see that if $\a=\l$ then $V$ is associative.
Thus $\a-\l\neq 0$. Furthermore, since $\phi(\U)$ is finitely generated over $F$ and finite dimensional over $K=F(\a)$, it is easy to see that $\a$ is algebraic over $F$ (see, for example, \cite[Lemma 1.1]{PS}).

\smallskip

If $\cha F=3$ and $\l=0$, then by Proposition \ref{prop5.1}  the superalgebra $\U(\a)$ has the ideal $I(\a)$ with $I(\a)^2=0$. Since $V$ is irreducible, $\phi(I)=0$, thus $\phi(\U)\cong \U(\a)/I(\a)=K\oplus K\,\bar a,\ \bar a^2=\a$.   It is clear that in this case $V$ as a $\phi(\U)$-module is isomorphic to the regular supermodule, which coincides with the bimodule $V_{\a,\l}^{\epsilon}(1|1)$ for $\epsilon=\a$ and $\l=0$.

\smallskip

Let $\cha F\neq 3$ or $\l\neq 0$; then the superalgebra $\U(\a)$ is simple and isomorphic to a generalized quaternion superalgebra over $K$,  with $\U(\a)_{\0}=K1+Kab,\ \U(\a)_{\1}=Ka+Kb$.  Therefore,
in this case $\phi(\U)=\U(\a)$, and we have the two possibilities, according to the structure of $\U(\a)_{\0}$.

If $\U(\a)_{\0}$ is a field, that is, the polynomial $f(u)=u^2+(\a-\l) u+\a(\a-2\l)$ is irreducible over $K$, then $V\cong\Reg\U(\a)$ as a right $\U(\a)$-supermodule, and  $V\cong V_{\a,\l}(2|2)$ as an $A$-bimodule.

If the polynomial $f(u)$ has a root $\varepsilon$ in $K$ then the $K$-subspace $L={\rm vect}\le ab-\epsilon,\ (\l-\a-\epsilon)\,a-\a\,b \re$ is a minimal graded right ideal of $\U(\a)$. It produces the irreducible superbimodule $V_{\a,\l}^{\epsilon}(1|1)$. Since $\tilde \U$ is simple and finite-dimensional, any other irreducible graded right $\U(\a)$-module up to the changing of parity is isomorphic to   $L$.
\ctd

\begin{Rem}\label{rem5.1}
Observe that  we have the following isomorphisms:
\bes
 V_{\a,\l}^{\epsilon_1}(1|1)^{\rm op}\cong V_{\a,\l}^{\epsilon_2}(1|1),\ \ \ V_{\a}(2|2)^{\rm op}\cong V_{\a}(2|2),
\ees
where $f(u)=(u-\epsilon_1)(u-\epsilon_2)$.
\end{Rem}
\begin{Rem}\label{rem5.2}
 If the field $F$ is algebraically closed  then we have $K=F$;  $\l\in \{0,\ 1\}$,
$f(u)$ reducible and $\dim_FV=2$. Moreover, if $\l=0$, we may take $\a=1$, and for $\a=0$ we may take $\l=1$. Thus up to the changing of parity, $V=V_{1}(1|1),\ V=V_{\a,1}^{\epsilon}(1|1)$, or $V=V_{1,0}^{\epsilon}(1|1)$, where in the last case $\epsilon^2+\epsilon+1=0$.

For $\l=1$ observe also that the  bimodules $V_1(1|1),\ V_{2,1}^0(1|1)$  are isomorphic to the bimodules of type 1,2 from Theorem \ref{th2.2}, and  $V_{\a,1}^{\epsilon}(1|1)$ is isomorphic to the bimodule of type 3 for the pair $(\a-1,\epsilon)$.
\end{Rem}

\begin{Rem}\label{rem5.3}
In case when $A=A_1=Fx,\ x^2=0$, we have the isomorphism $\U(A)\cong \U_1(A^{\sharp})$, hence $A$ has the same irreducible superbimodules as $A^{\sharp}$ does.
\end{Rem}

\subsection{$A_{\1}^2=0$, the general case.}

We will consider now the general case of superalgebras with $A_1^2=0$.

\begin{Th}\label{th5.2}
Let $V$ be an irreducible faithful alternative bimodule over the superalgebra $A$ with $A_{\1}^2=0$. Then there exist a field  extension $K$ of $F$, a nonzero algebraic element $\a\in K$, an $F$-subspace $T$ and an $F$-subalgebra $S$ of $K$   such that:
\begin{itemize}
\item[$i)$] $V$ is a unital  irredicuble faithful alternative bimodule over the superalgebra $B=K\oplus Ku,\ u^2=0$, and the structure of $V$ as a $B$-bimodule is given by Theorem \ref{th5.1}, with $\r(u)^2=\a$.
\item[$ii)$] $S+Tu$ is an $F$-subsuperalgebra of $B$ which is isomorphic to $A$.
\item[$iii)$]   $K=\alg_F\le \a T^2\re$.
\item[$iv)$] The action of $A$ on $V$ is inherited by that of $B$.
\end{itemize}
Conversely, let $K$ be a field extension of $F$, and let  $T,\,S$ and $\a$  be an $F$-subspace,  $F$-subalgebra and $F$-algebraic element of $K$ such that $ST+TS\subseteq T$ and $\alg_F\le \a T^2\re=K$. Consider the superalgebra $B=K\oplus Ku,\ u^2=0$, then $A=S+Tu$ is a subsuperalgebra of $B$ and every irreducible faithful alternative birepresentation $(\l,\r)$ of $B$ with $\r(u)^2=\a$ is so over $A$.
\end{Th}
\Proof
Consider the subsuperalgebra $\G$ of the centroid ${\mathcal E}={\mathcal E}(E)$ of the split null extension $E=V\oplus A$: $\G=\{\a\in {\mathcal E}\,|\,V\a\subseteq V,\ A\a\subseteq A\}$. By Proposition \ref{prop5.1}, $\G$ is supercommutative and isomorphic to a subsuperalgebra of the centralizer $\ZZ=\ZZ(V)$ which is a division superalgebra. We will identify $\G$ with its image in $\ZZ$.

Let $(L,R):A\rightarrow\End E$ be the regular birepresentation of $A$ in $E$,
then for any $x,y\in A_{\1}$ the operator $D=R_xR_y$ is a derivation of the superalgebra $E$ due to (\ref{id2.2}). Since $D|_A=0$,
it follows that $D\in\G\subseteq\ZZ$.

If $R_xR_y=0$ for all $x,y\in A_{\1}$ then $(VA_{\1})A_{\1}=0$.
But $VA_{\1}$ is a subbimodule of $V$, and if $V=VA_{\1}$ then
$V=(VA_{\1})A_{\1}=0$, a contradiction. Thus $VA_{\1}=0$.
It follows from (\ref{Ialt}) that $L_xL_y=R_xR_y=0$, hence we have
$A_{\1}(A_{\1}V)=0$ and $A_{\1}V=0$ as well. Thus $A_{\1}\subseteq \Ann V=0$, a contradiction.

Let $R_xR_y\neq 0$, show that then  $R_x^2\neq 0$. In fact, we have seen that $0\neq D=R_xR_y\in \ZZ$. Therefore, $V=VD=VD^2=VR_xR_yR_xR_y$. But $R_xR_y=R_yR_x$ by (\ref{Ialt}), hence
$V=VR_x^2R_y^2$ and $R_x^2\neq 0$. Since $R_x^2$ is invertible on $V$, so is $R_x$.

Assume now that $R_zR_t=0$ for some $z,t\in A_{\1}$, then $R_t$ is not invertible and hence $R_uR_t=L_uL_t=0$ for any $u\in A_{\1}$.
Thus $VA_{\1}R_t=0=(A_{\1}V)L_t$ and $t=0$. We proved that for any $0\neq x\in A_{\1}$ the operator $R_x$ is invertible on $V$.  In particular, we have $V_{\1}=V_{\0}x,\ V_{\0}=V_{\1}x$.

Let below $v\in V,\ x,y\in A_{\1},\ a,b,c\in A_{\0}$.
From super-linearized identity (\ref{id2.9}) we have
\bes
(vx)(ay)\pm (yx)(av)=(v(xa)y\pm (y(xa))v,
\ees
which gives $(vx)(ay)=(v(xa))y$ or
\bee\label{id5.1}
R_xR_{ay}=R_{xa}R_y.
\eee
For $y=x$ we have $R_xR_{ax}=R_xR_{xa}$, which by invertibility of $R_x$ implies $R_{xa}=R_{ax}$ or $R_{[a,x]}=0$, and eventually $[a,x]=0$, that is,
\bee\label{id5.2}
[A_{\0},A_{\1}]=0.
\eee
Furthermore, applying (\ref{id5.1}) and (\ref{id5.2}), we get
\bes
R_{(ax)b}R_x=R_{ax}R_{bx}=R_{ax}R_{xb}=R_{xb}R_{ax}=R_{(xb)a}R_{x}=R_{a(xb)}R_{x},
\ees
which implies $R_{(a,x,b)}=0$ and therefore
\bee\label{id5.3}
(A_{\0},A_{\1},A_{\0})=0.
\eee

Let $R$ be the $F$-subalgebra of $\ZZ_0$ generated the elements $\{R_xR_y\,|\,x,y\in A_{\1} \}$, then $R\subseteq \G_{\0}$ is a commutative domain. Denote by $K$ the field of fractions of $R$.
Fix $0\neq x\in A_{\1}$ and let $\a=(R_x)^2$. Then for any $y\in A_{\1}$, $R_y= (R_xR_y)\a^{-1}R_x$. Consider the $F$-linear mapping $\tau:A_{\1}\rightarrow K$ given by $\tau(y)=(R_xR_y)\a^{-1}$, then $R_y=\tau(y)R_x$. Clearly, $\ker \tau=0$, hence $T=\tau(A_{\1})$ is an $F$-vector subspace of $K$ isomorphic to $A_{\1}$.  Furthermore, define an $F$-linear mapping $\s:A_{\0}\rightarrow K$ by $\s(a)=\tau(ax)$ and denote $S=\s(A_{\0})$.

\begin{Lem}\label{lem5.1}
The mappings $\tau$ and $\s$ satisfy the conditions:
\bes
 i)& \tau(ay)=\s(a)\tau(y),  &ii)\, \s(ab)=\s(a)\s(b),\\
 iii)& \ker\s=0, &iv)\, L_y=\tau(y)L_x,
\ees
for all $a,b\in A_{\0},\ y\in A_{\1}$.
\end{Lem}
\Proof We have
\bes
\s(a)\tau(y)&=&(R_xR_{ax}\a^{-1})(R_xR_{y}\a^{-1})=R_{ax}R_{y}R_xR_{x}\a^{-2}\\
&=&R_{ax}R_{y}\a^{-1}\stackrel{(\ref{id5.1})}=R_{x}R_{ay}\a^{-1}=\tau(ay),
\ees
proving $i)$. Furthermore,
\bes
\s(a)\s(b)&=&\tau(ax)\tau(bx)=(R_xR_{ax}\a^{-1})(R_xR_{bx}\a^{-1})=R_{ax}R_{bx}\a^{-1}\\
&\stackrel{(\ref{id5.1})}=&R_xR_{a(bx)}\a^{-1}\stackrel{(\ref{id5.3})}=R_xR_{(ab)x)}\a^{-1}=\s(ab),
\ees
which proves $ii)$. Let now $a\in\ker\s$, then by $i)$, $\tau(ay)=0$ for any $y\in A_{\1}$ and $aA_{\1}=0$. By Lemma \ref{lem3.1}, $(V,a,x)$ is a subbimodule of $V$.
If $(V,a,x)=0$ then we have $R_aR_x=0$ and $R_a=0$. Similarly, $L_a=0$ and  $a\in \Ann V=0$. If $(V,a,x)=V$ then $Va=(V,a,ax)=0$ and $aV=0$, and again $a=0$.

To prove $iv)$, consider
\bes
L_y&=&L_y\a\a^{-1}=L_yR_xR_x\a^{-1}=L_yL_xL_x\a^{-1}\\
&=&R_yR_xL_x\a^{-1}=\tau(y)L_x;
\ees
\ctd
\begin{Cor}\label{cor5.1}
Let $B=K\oplus Ku,\ u^2=0$, with $B_{\0}=K,\ B_{\1}=Ku,$ then the
 application $a_0+y_1\mapsto \s(a_0)+\tau(y_1)u$ is an isomorphism    of $A$ and $S+Tu\subseteq B$. In particular, $A_{\0}$ is an associative and commutative domain.
\end{Cor}

\begin{Lem}\label{lem5.2}
For any $a\in A_{\0},\ s\in A_{\0}\cup A_{\1}$ hold
\bes
i)\, R_a=L_a\in \ZZ, \ \ ii)\, R_aR_s=R_{as}.
\ees
\end{Lem}
\Proof
 By Lemma \ref{lem3.1},
$(V,a,b)$ is a subbimodule of $V$. If it is not zero then by (\ref{id2.4})  we have $V=(V,a,b)=((V,a,b),a,b)=[a,b](V,a,b)=0$, a contradiction. Therefore $R_aR_b=R_{ij}=R_{ba}=R_bR_a$.

Furthermore, in view of (\ref{id2.2}), (\ref{id5.2}) and (\ref{id5.3}), $(V,a,x)$ is a subbimodule of $V$ for any $x\in A_{\1}$.
If $(V,a,x)=0$ then $R_aR_x=R_{ax}=R_{xa}=R_xR_a$. Moreover, in this case we have
$\s(a)R_x=R_{ax}=R_aR_x$, implying $R_a=\s(a)$. Similarly $L_a=\s(a)$ and $L_a=R_a\in\ZZ$.

Assume that $(V,a,x)=V$.
Then for any $v\in V$ we have $(v,a,x)a=(v,a,ax)=(v,a,xa)=a(v,a,x)$, that is, $R_a=L_a$. Furthermore, we have
 $R_a\circ R_x=2R_{ax}$, $[R_a,R_x]=[L_a,R_x]=-R_aR_x+R_{ax}$.
Summing the two equations we get $3R_aR_x=3R_{ax}$. By (\ref{id2.2}), $R_aR_x-R_{ax}\in \G_{\1}\subseteq \ZZ_{\1}$.  Since $\G$ is a supercommutative domain, $\G_{\1}$ may be non-zero only if $\cha F=2$.
Hence we have $R_aR_x=R_{ax}$.
\ctd
\begin{Lem}\label{lem5.3}
If $cha F\neq 3$ then $L_xR_x\not\in\ZZ$.  If $L_xR_x\in\ZZ$ then $\cha F=3$ and $R_x=L_x$.
\end{Lem}
\Proof
Assume that $L_xR_x\in\ZZ$, then we have $L_xR_xR_x=R_xL_xR_x$ which implies $L_xR_x=R_xL_x$.
By (\ref{Ialt}), $2L_xR_x=[L_x,R_x]_s=-(R_x)^2=-(L_x)^2$. Therefore,
$2L_x=-R_x,\ 2R_x=-L_x$, implying $L_x=R_x$. Since $x\in A_{\1}$ is arbitrary, $V$ is supercommutative, and we have
\bes
3vR_xR_y=3(v,x,y)=3(x,y,v)\stackrel{(\ref{id2.1})}=[xy,v]-x[y,v]-(-1)^{|y||v|}[x,v]y=0,
\ees
which is impossible  if $\cha F\neq 3$. If $\cha F=3$ then we have $(L_xR_x-\a)^2=-\a^2-\a L_xR_x-2\a L_xR_x+\a^2=0$ (see Section 5.1). Therefore, if $L_xR_x\in\ZZ$ then $L_xR_x=\a=R_xR_x$, and $R_x=L_x$.
\ctd

\smallskip
Return to the proof of the Theorem. Assume that $L_xR_x\not\in \ZZ$ then there exists $v\in V_{\0}$ such that $v$ and $w=vL_xR_x$ are linearly independent over the division algebra $\ZZ_{\0}$. For $y,z\in A_{\1}$ we have $L_yR_z=(R_yR_z\a^{-1})(L_xR_x)\in R[\a^{-1}]L_xR_x$.
Since $V$ is irreducuble, it follows from Lemma \ref{lem5.2} and (\ref{Ialt})  that $vR+wR[\a^{-1}]=V_{\0}$.  Thus we have
\bes
V_{\0}=vR+wR[\a^{-1}]\subseteq vK+wK\subseteq v\ZZ_{\0}+w\ZZ_{\0}\subseteq V_{\0},
\ees
which implies that $R=K=\ZZ_{\0}$.

If $L_xR_x\in\ZZ$ then $\cha F=3,\ L_x=R_x$, and for any $0\neq v\in V_{\0}$ we have $V_{\0}=vR$, which implies as above  $R=K=\ZZ_{\0}$. Notice that $R_yR_z=\tau(y)\tau(z)\a\in T^2\a$, hence $K=\alg\le T^2\a\re$. 

\smallskip

Define a unital $K$-linear action of $B$ on $V$ by setting
$v\cdot u=vR_x,\ u\cdot v=vL_x$, then evidently $V$ is an alternative irreducible faithful $B$-bimodule. The restriction of this action on $S+Tu\cong A$ coincides with the $A$-bimodule structure of  $V$.

\smallskip

Conversely, let $K$ be a field extension of $F$, and let  $T,\,S$ and $\a$  be an $F$-subspace,  $F$-subalgebra and $F$-algebraic element of $K$ such that $ST+TS\subseteq T$ and $\alg_F\le \a T^2\re=K$. Consider the superalgebra $B=K\oplus Ku,\ u^2=0$, then $A=S+Tu$ is a subsuperalgebra of $B$. Let $V$ be an irreducible unital faithful alternative $B$-superbimodule such that $\r(u)^2=\a$, and assume that $W\subseteq V$ is an $A$-subsuperbimodule of $V$. We have for any $w\in W,\ t_1,t_2\in T$: $(w\cdot t_1u)\cdot t_2u=(t_1t_2) (w\cdot u)\cdot u =(\a t_1t_2)w \in W$, hence  $W$ is a $K$ subspace of $V$.  Let $0\neq t\in T$, then $W\supseteq t^{-1} W\cdot (tu)=W\cdot u$, and $W=V$.
\ctd

\section{Bimodules over the simple superalgebras.}
\subsection{Reduction of prime  superalgebras to simple ones}

\hspace{\parindent}
We start this section with reduction of the case of prime associative superalgebras to the simple finite dimensional ones.

We will say that a superalgebra $A$ with the center $Z(A)=Z_{\0}\oplus Z_{\1}$ is {\it an even central order} in a superalgebra $B$ if $B=(Z_{\0}^*)^{-1}A$. In this case we will also call $B$ to be {\it the even central closure} of $A$.

\begin{Lem}\label{lem6.1}
Let $A$ be a  prime associative superalgebra and let $V$ be a faithful irreducible alternative $A$-superbimodule which is not associative. Then $Z_{\0}=Z(A)_{\0}\neq 0$, and the even central closure $B=(Z_{\0}^*)^{-1}A$ is a simple superalgebra which is  finite dimensional over the even part of its supercenter $(Z_{\0}^*)^{-1}Z$.
Moreover, $(Z_{\0}^*)^{-1}Z\subseteq \ZZ(V)$, and $V$ has a natural structure of a $B$-bimodule.
\end{Lem}
\Proof Since $A$ is associative,  we can repeat the arguments from Section 4 and to get $V_{as}:=\{v\in V\,|\,(v,A,A)=0\}=0$ and $V\cdot [N,A_{\0}]=0$, where $N$ is the associative center of the split null extension $E=A\oplus V$. Therefore $A_{\0}$ satisfies the identity $[[x,y]^4,z]=0$. By \cite{BGR}, $A$ is also a $PI$-algebra.

Consider the quotient algebra $A/\cal{B}(A),$ where $\cal{B}(A)$ is the Baer radical of the algebra $A$. Since $A/\cal{B}(A)$ is semiprime (as algebra) and $PI$, it follows that there exists a nontrivial central polynomial $g(x_1,\ldots,x_n)$ of $A/\cal{B}(A)$ \cite{Rowen}. Hence $\left[g(a_1,\ldots,a_n),A\right]\subseteq \cal{B}(A)$ for arbitrary homogeneous elements $a_k.$ But the Baer radical of a prime superalgebra does not contain nonzero homogeneous elements (see \cite{CohenMontgomery}), therefore $\left[g(a_1,\ldots,a_n),A\right]=0$ and $Z_{\0}\neq 0.$

Consider the prime superalgebra $B=(Z_{\0}^*)^{-1}A.$ Let $I$ be a nonzero graded ideal of $B.$ Then $I$ is a prime superalgebra as well, and $I\bigcap Z_{\0}(B)=Z_{\0}(I)\neq 0.$ The proof of these statements repeats verbatim the non-graded case.  Hence $I=B$ and $B$ is a simple superalgebra. It is also clear that $Z_{\0}(B)=(Z_{\0}^*)^{-1}Z$ is a field.

Due to \cite{Wall},  $B$ is  simple as an algebra or there exists a proper (non-graded) ideal $I$ in $B$ such that  $B=\pi_0(I)\oplus\pi_1(I),$ where $\pi_i:I\longrightarrow B_i$ are the natural projections. It is easy to see that in this case  $B\cong B_{\0}[u]$ where $B_{\0}$ is simple.

In the first case, since $B$ is a simple PI-algebra it follows that $B$ is finite-dimensional over $Z(B).$ But $Z(B)$ is finite-dimensional over $Z_{\0}(B).$  Indeed, if $Z_{\1}(B)\neq 0$ then, fixing $0\neq z_1\in Z_{\1}(B),$ we have $z_1^2=\alpha\neq 0\in Z_{\0}(B)$ and then for any $x\in Z_{\1}(B)$ we get $x=x(\alpha^{-1}z_1^2)=(\alpha^{-1}xz_1)z_1\in Z_{\0}(B)\cdot z_1.$

In the second case $B_{\0}$ is a simple $PI$-algebra, hence $B_{\0}$ is finite-dimensional over $Z(B_{\0}).$ It is easily seen that $B$ is finite-dimensional over $Z(B_0)$ too and, in this case, $Z(B)=Z(B_{\0})[u].$

\smallskip

By Lemma \ref{lem4.1}, $Z(A)$ is contained in the centralizer $\ZZ(V)$ of the bimodule $V$, which is a graded division superalgebra. Therefore, $(Z_{\0}^*)^{-1}Z\subseteq \ZZ$ and $V$ has a natural structure of a $B$-bimodule. Clearly, $V$ is irreducible and faithful over $B$.
\ctd

\smallskip

In non-associative case we have

\begin{Lem}\label{lem6.2}
Let $A$ be a  prime alternative non-associative non-trivial superalgebra of $\cha\neq 3$ and let $V$ be a faithful irreducible alternative $A$-superbimodule. Then
$\cha F=2$,  $Z_{\0}=Z(A)_{\0}\neq 0$, and the even central closure $B=(Z_{\0}^*)^{-1}A$ is a simple superalgebra of type $\OO(4|4)$ or $\OO[u]$.
As above, $(Z_{\0}^*)^{-1}Z\subseteq \ZZ(V)$, and $V$ has a natural structure of a $B$-bimodule.
\end{Lem}
\Proof The first part of the Lemma follows from the classification of prime alternative superalgebras in \cite{Sh2}. The second part follows from Lemma \ref{lem4.1} as above.
\ctd

Therefore, except for the case of $\cha 3$, the classification of irreducible
alternative non-associative superbimodules over prime superalgebras  is reduced to bimodules over simple central finite-dimensional superalgebras, where we call a superalgebra $A$ over a field $F$  {\em central} if $Z(A)_{\0}=F$.


\smallskip

We will show next that, modulo the superalgebras of type $B=K[u]$, for the case of $\cha \neq 2$  it suffices to consider only non-associative superalgebras.
\begin{Lem}\label{lem6.3}
In the conditions of Lemma \ref{lem6.1}, if  $\cha F\neq 2$ then the even central closure $B$ of $A$ is of the type $B=K+Ks,\ s^2=\l\in K,\ \l\neq 0$, and the structure of $V$ and $A$ is given by Theorem \ref{th5.1}.
 \end{Lem}		
\Proof By Lemma \ref{lem6.1}, $Z=Z_{\0}(A)\neq 0$, $A$ is an even central order in a simple central superalgebra $B$ over the field $K=(Z_{\0}^*)^{-1}Z_{\0}$, and $V$ is a faithful irreducible alternative bimodule over $B$. Let $\tilde K$ be the algebraic closure of $K$; then $\tilde B=\tilde K\otimes B$ is a simple central superalgebra over $\tilde K$ and $\tilde V=\tilde K\otimes V$ is an alternative unital bimodule over $\tilde B$. If $\dim_{\tilde K}{\tilde B}>2$, then by \cite{Pisar} $\tilde V$ is completely reducible and  all its irreducible components are associative. Hence in this  case $V$ is associative, contrary to the assumption. Therefore, $\dim_{\tilde K}{\tilde B}=\dim_{K}{B}=2$ and $B=K+Ks,\ s^2=\l\in K^*$.
\ctd

\begin{Rem}\label{rem6.1}
Observe that, contrary to the case of algebras in Section $4$, $A$ may be not simple and hence in general $B\neq A$. As an example, consider the superalgebra $B=\QQ[\sqrt{1}]=\QQ\oplus\QQ u,\ u^2=1$ and let $V=\QQ v_0\oplus\QQ v_1$ be the irreducible supermodule of type 3 from Theorem  \ref{th2.2} for $\a=-\tfrac37,\ \d=\tfrac87$. Let  $A=\QQ_7[\sqrt{1}]=\QQ_7\oplus\QQ_7u,\ u^2=1$, where $\QQ_7$ stands for the subring of rational numbers whose denominators are not divisible by 7. Then $V$ is an irreducible faithful alternative $A$-bimodule, and $B=(Z_{\0}^*)^{-1}A\neq A$.
\end{Rem}
If one wanted $A$ to be an algebra over a field and not just a ring, then one may take in this example $A=F[[x]][\sqrt{1}],\ B=F((x))[\sqrt{1}],\ \a=\tfrac{-x^2+1}{x^2+x+1},\ \d=x(1+\a)$  where $F$ is a field, $F[[x]]$ and $F((x))$ are the algebra of formal series and the field of formal Laurent series over $F$.

\subsection{Superbimodules in  characteristic 2.}

 \hspace{\parindent}
In this section, we will classify irreducible alternative non-associative superbimodules of characteristic 2.

Recall that every $\Z_2$-graded alternative algebra of characteristic 2 is an alternative superalgebra. The  simple superalgebras $\OO(4|4)$ and $\OO[u]$ provide us the irreducible  superbimodules $\Reg\OO(4|4),\ \Reg(\OO[u])$.

Furthermore,  consider a generalized quaternion algebra $\HH$
 with the $\Z_2$-grading coming from the Cayley-Dickson process; then it is an alternative superalgebra which we denote as $\HH(2|2)$.
As in the non-graded case, the category of $\Z_2$-graded Cayley bimodules over $\HH(2|2)$ is isomorphic to the category of (graded) left $\HH(2|2)$-modules. For a left  $\HH(2|2)$-module $L$ we again denote by $\Cay L$ the corresponding Cayley $\HH(2|2)$-bimodule. Then $\Cay L$  is an alternative superbimodule over $\HH(2|2)$ in the characteristic 2 case.

Finally, let $A$ be an alternative algebra and $V$ be an $A$-bimodule; consider the split null extension $E=E(A,V)$. The double  $E[u]=E+Eu,\ u^2=\a\neq 0$ is an alternative superalgebra with the even part $E$ and the odd part $Eu$ (recall that $\cha F=2$). Moreover, the ideal $V[u]=V+V u$ may be considered as a bimodule over the alternative superalgebra $A[u]=A+A u$, and if $V$ is irreducible then so is $V[u]$. We will denote this bimodule as $V[u]$ or as $V[\sqrt1]$ when $u^2=1$. 


\begin{Th}\label{th6.1}
Let $V$ be an irreducible faithful alternative non-associative bimodule over a nontrivial prime alternative superalgebra $A$  of characteristic $2$ with $\dim_{(Z^*)^{-1}Z}(Z^*)^{-1}A>2$. Then, up to the graded isomorphism, we have the following possibilities:
\begin{itemize}
\item $A=\HH(2|2),\ V=\Cay L$,  $L$ is a minimal $\Z_2$-graded left ideal of $\HH(2|2)$;
\item $A=\HH[u],\ V=(\Cay L)[u]$,  $L$ is a minimal left ideal of $\HH;$
\item $A=\OO(4|4),\ V=\Reg (\OO(4|4))$;
\item $A=\OO[u],\ V=\Reg (\OO[u])$.
\end{itemize}
\end{Th}
\Proof
Consider first the associative case. By Lemma \ref{lem6.1},  $B=(Z^*)^{-1}A$ is a simple central superalgebra over the field $K=(Z^*)^{-1}Z$. Let $\tilde K$  be the algebraic closure of $K$, then $\tilde B=\tilde K\otimes_K B$ is a simple central superalgebra over $\tilde K$ and by \cite{Wall} $\tilde B$ is isomorphic to $M_{m|k}(\tilde K),\ 1\leq k\leq m,$ or to $M_n(\tilde K)[\sqrt1],\ n>1$.
Let us show that it should necessary hold $m+k=n=2$. Consider $\tilde V=\tilde K\otimes_K V$; it is a bimodule over $\tilde B$ and by \cite[Lemma 7.5]{ZSSS}  we have $[N,\tilde B]_s\subseteq \Ann (\tilde V,\tilde B,\tilde B)$, where $N$ is the associative center of the split null extension $E(\tilde B,\tilde V)$. Notice that since $B$ is associative and $V$ is irreducuble, we have $V=(V,B,B)$. Therefore $V[N,\tilde B]_s=0$, which implies easily that $[N,\tilde B]_s=0$. In particular, the following inclusion holds in $\tilde  B$:
\bee\label{id6.1}
[[a,x]\circ_s[a,y],a]\in Z_{\0}(\tilde B)=\tilde K
\eee
for any $a\in \tilde B_{\0},\ x,y\in\tilde B_{\1}$. In fact, by \cite[(7.25) and Lemma 7.1]{ZSSS}, we have the inclusion $[[a,b]^2,a]\in N(A)$
in any alternative algebra $A$. Superlinearizing this inclusion on $b$, we get  the inclusion
\bee\label{id6.2}
[[a,x]\circ_s[a,y],a]\in N(E),
\eee
for any $a\in E_{\0},\ x,y\in E_{\0}\cup E_{\1}$, which implies (\ref{id6.1}). Assume now that $\tilde B=M_{m|k}(\tilde K)$ with $m>1$ and substitute in (\ref{id6.1}) $a=e_{11}+e_{22}+e_{12},\ x=e_{1,m+1},\ y=e_{m+1,1}$. We obtain
\bes
[[a,x]\circ_s[a,y],a]&=&[[e_{1,m+1},e_{m+1,1}+e_{m+1,2}],a]\\
&=&[e_{11}+e_{12}+e_{m+1,m+1},e_{11}+e_{22}+e_{12}]=e_{12}\not\in\tilde K,
\ees
a contradiction. Similarly, assume that $\tilde B=M_n(\tilde K)[\sqrt1]$ with $n>2$ and substitute in (\ref{id6.1}) $a=e_{11}+e_{22}+e_{12},\ x=e_{13}u,\ y=e_{31}u$:
\bes
[[a,x]\circ_s[a,y],a]&=&[[ue_{13},u(e_{31}+e_{32})],a]\\
&=&[e_{11}+e_{12}+e_{33},e_{11}+e_{22}+e_{12}]=e_{12}\not\in\tilde K,
\ees
a contradiction again.

Returning to the superalgebra $B$ we conclude that in the first case $B$ is a generalized quaternion superalgebra $\HH(2|2)$ and in the second case $B_{\0}\cong\HH$.

\smallskip

Let $B=\HH(2|2),\ B_{\0}=K+Ki,\ B_{\1}=Kj+Kk,\ i^2=i+\a,\ k=ij=j(i+1),\ j^2=\b,$
where $\a,\b\in K,\ (4\a+1)\b\neq 0$ (see \cite[Chapter 2]{ZSSS}). Substituting in (\ref{id6.2}) $a=i,\ x=j,\ y=v\in V$, we get
\bes
([[i,j]\circ_s [i,v],i],B,B)=([[[v,i],j],i],B,B)=0,
\ees
that is,  $[[[V,i],j],i]\subseteq V_{as}=\{v\in V\,|\, (v,B,B)=0\}$ which is a subbimodule of $V$. Since $V$ is nonassociative, $V_{as}=0$, and we have
$[[[V,i],j],i]=0$. Let $W=\{w\in V\,|\,[w,i]=0\}$, show that $W=0$. For any $w\in W,\ b\in B$ we have
\bes
(w,b,i)=(w,b,i^2-\a)=(w,b,i^2)=(w\circ i,b,i)=([w,i],b,i)=0,
\ees
that is, $(W,B,i)=0$. Further, $(W,B,B)\circ i\subseteq (W\circ i,B,B)+W\circ (i,B,B)=0$, thus $(W,B,B\circ i)\subseteq (W,B,B)\circ i+(W,B,i)\circ B=0$.
Observe that $j=j\circ i,\ k=k\circ i$, which implies that $(W,B,B)=0$ and $W\subseteq V_{as}=0$.

Therefore, we have $[[V,i],j]=0$ and similarly $[[V,i],k]=0$. Let us prove that $U=[V,i]=V\circ i$ is a subbimodule of $V$. Clearly, $Ui+iU\subseteq U$. Note that $(V,B,i)=(V,B,i^2)= (V,B,i)\circ i\subseteq U$. We have also $[U,j]=[U,k]=0$. Consider
\bes
[v,i]j&=&[v,i^2]j=[v\circ i,i]j=[v\circ i,ij]-i[v\circ i,j]+3(v\circ i,i,j)\\
&=&[v\circ i,k]+(v\circ i,i,j)=(v\circ i,i,j) \in (V,B,i)\subseteq U.
\ees
Thus $Uj\subseteq U$ and similarly $Uk\subseteq U$, hence $U$ is a subbimodule of $V$. We have already showed that the mapping $v\mapsto v\circ i$ is injective. Therefore $U=V$, and we have $[V,j]=[V,k]=0$. Consider
$3(v,j,j)=[vj,j]-v[j,j]+[v,j]j=0$. Similarly, we have $(V,k,k)=0$, which proves that $V$ is an alternative $\Z_2$-graded $B$-bimodule.

Denote $T=\{v-[v,i]\,|\,v\in V\}$, then clearly $Ti+iT\subseteq T$. Furthermore,
\bes
[v,i]j&=&[v,ij]-i[v,j]+3(v,i,j)=3(v,i,j)\\
&=&[vi,j]-v[i,j]-[v,j]i=v[i,j]=vj,
\ees
and similarly $[v,i]k=vk$. Thus $Tj=Tk=0$ and $T$ is a subbimodule of $V$. Since $V$ is faithful, $T=0$ and we have $v\circ i=v$ for any $v\in V$. Therefore, $bv=v\bar b$ for any $b\in B$, and $V$ is a Cayley bimodule over $B$. Since $V$ is irreducible, $V=\Cay L$, where $L$ is an irreducible left $\HH(2|2)$-module. It is isomorphic to a minimal graded left ideal of $B=\HH(2|2)$. As at the end of the proof of Theorem \ref{th4.1}, we get that $L$ is an irreducible graded left $A$-module and $A$ is a graded primitive $PI$-algebra. By the graded version of Kaplansky's theorem, $A$ is graded simple. Hence $Z(A)$ is a graded division algebra, $K=Z_{\0}$ and $A=B=\HH(2|2)$.

\smallskip

Consider now the case when $\tilde B=M_2(\tilde K)[\sqrt1]=M_2(\tilde K)+M_2(\tilde K)u$. We have $Z_{\1}(\tilde B)=\tilde Ku\neq 0$, therefore  $Z_{\1}(B)\neq 0$ and $B$ is not simple as an algebra. Consequently $B=B_{\0}[u]=\HH[u]=\HH+\HH u, u^2=\a\neq 0$. Take $a\in \HH$ and consider $W=(V,a,u)$. It follows from (\ref{id2.2}) that $W$ is a subbimodule of $V$. Assume that $W=V$, then $V=((V,a,u),a,u)$. But we have for any $v\in V$ by super-linearized (\ref{id2.4})
\bes
((v,a,u),a,u)=((u,a,u),a,v)+[a,u](v,a,u)+[a,v](u,a,u)=0,
\ees
a contradiction. Therefore $(V,\HH,u)=0$, and by (\ref{id2.3})
$[u,(V,\HH,\HH)]\subseteq (u,V,\HH)=0$.
Let us show that if $[u,v]=0$ then $(u,u,v)=0$. In fact, by (\ref{id2.1}),
$3(u,u,v)=[u^2,v]-u[u,v]-[u,v]u=0$.  Thus $((V,\HH,\HH),u,u)=0$.

Consider now the space $(V,u,u)$; it is also a subbimodule of $V$. Assume that $(V,u,u)=V$, then we have
$$
(V,\HH,\HH)=((V,u,u),\HH,\HH)=((V,\HH,\HH),u,u)=0.
$$
We have already seen that $(V,\HH,u)=0$ as well. Take $a,b\in \HH$ and consider
\bes
(v,a,bu)\stackrel{(\ref{id2.7})}=(v,u,ba)+u(v,a,b)+a(v,u,b)=0.
\ees
Therefore $(V,\HH,B)=0$, which imlies by (\ref{id2.2}) that $[V,\HH]\subseteq V_{as}=0$. Consequently, $V$ is an associative and commutative $\HH$-bimodule. But it is impossible since for any $a,b\in \HH$ and $v\in V$ we would have
$$
(ab)v=a(bv)=a(vb)=(vb)a=v(ba)=(ba)v,\  [a,b]v=v[a,b]=0,
$$
and $[a,b]=0$, a contradiction. Therefore, $(V,u,u)=0$.

It easy now to see  that $(V,u,B)=0$ and then $[V,u]\subseteq V_{as}=0$.
Therefore, $u\in Z(E(V,B))$ and $V=V_{\0}+V_{\0}u=V_{\0}[u]$.
Clearly, $V_{\0}$ is an irreducible alternative $\HH$-bimodule. If $V_{\0}$ were associative then $V$ would be associative as well, which is not the case. Therefore, $V\cong \Cay L$ where $L$ is a minimal left ideal of $\HH$. As in the previous case, $A_{\0}$ is a primitive $PI$-algebra and $Z(A_{\0})=Z(A)_{\0}$ is a field. Therefore $A=B$.

\smallskip

Assume now that $B=\OO[u]=\OO+\OO u,\ u^2=\l\in K,\ \l\neq 0$, where $\OO$ is a central Cayley-Dickson algebra over $K$.
As above, we have $(V,u,\OO)=[u,(V,\OO,\OO)]=((V,\OO,\OO),u,u)=0$. Assume that $(V,u,u)=V$, then $(V,\OO,\OO)=0$ and $V$ is a unital associative bimodule over $\OO$, which is impossible. Hence $(V,u,u)=0$. As above, we have $(V,u,\OO)=[V,u]=0$, hence $u\in Z_{\1}(E(V,B))$ and $V=V_{\0}[u]$. It is clear now that $V_{\0}\cong\Reg\OO$. As in the proof of Theorem \ref{th4.1}, we have also $A=B$.

\smallskip

Let finally  $B=\OO(4|4)$.
Then $B_{\0}=\HH$ and $B_{\1}$ is a Cayley bimodule over $\HH$. By \cite[Lemma 12]{Sh2}, $V=V_a\oplus V_c$, where $V_a$ is an associative $\HH$-bimodule and $V_c$ is a Cayley $\HH$-bimodule. Moreover, by \cite[Lemma 3.2]{ConchShest}, we have
$V_aB_{\1}+B_{\1}V_a\subseteq V_c,\  V_cB_{\1}+B_{\1}V_c\subseteq V_a$. Since $V$ is faithful, this implies that $V_a\neq 0,\ V_c\neq 0$. Let $Z_a=\{v\in V_a\,|\,[v,\HH]=0\}$. As in the proof of Lemma 3.3 in \cite{ConchShest} we have that $Z_a\neq 0$ and $[Z_a,B]=0$. Furthermore, $3(Z_a,B,B)\subseteq [BB,Z_a]+B[B,Z_a]+[B,Z_a]B=0.$
Hence $0\neq Z_a\subseteq Z(E(V,B))$. Choose some homogeneous element $0\neq u\in Z_a$, then the  subspace $u\cdot B$ is a $B$-subbimodule of $V$ and the mapping $\phi : a\mapsto u\cdot a$ is a $B$-bimodule homomorphism of $\Reg B$ onto $uB$,
in case when $u$ is even, or of $(\Reg B)^{op}$ onto $uB$, in case when $u$ is odd. Since both $\Reg B$ and $(\Reg B)^{op}$ are irreducible, and $\phi(1) = u\neq 0$, we have that $uB = V$ is isomorphic to $\Reg B$ or to $(\Reg B)^{op}$. Clearly, in both cases we have $A=B$.
\ctd

\subsection{Superbimodules over the  superalgebra $B(\Gamma,D,\gamma)$.}

\hspace{\parindent}

In this section we consider the last remained case of superbimodules in characteristic 3. Unlike the other cases, we can not classify  irreducible alternative superbimodules of characteristic 3 in any dimension, since the prime alternative superalgebras of characteristic 3 are classified only under certain restrictions.

Nevertheless, due to \cite{ConchitaShest}, every finite dimensional non-associative   prime alternative superalgebra $A$ is simple.
Thus, if $\cha A=3$ then $A$ is one of the superalgebras $B(1|2),\ B(4|2),\ B(\G,D,\g)$. The supermodules over the superalgebras
$B(1|2),\ B(4|2)$ were classified in \cite{ConchShest, T4} (see Theorems \ref{th2.3}, \ref{th2.4}). Therefore, to finish the classification of irreducible finite dimensional superbimodules,
 it suffices to classify them over the superalgebra $B(\G,D,\g)$. We will do this for any dimension.

\smallskip

Below $V$ denote an irreducible faithful bimodule over the  superalgebra $B=B(\G,D,\gamma)$ which is not associative.

We start with the studying the structure of $V$ as a bimodule over $B_{\0}=\G$.

\begin{Lem}\label{lem6.4}
$(V,\G,\G)=[V,\G]=0.$
\end{Lem}
\Proof Fix $a,b\in \G.$ Since $\G$ is commutative and $(B,\G,\G)=0,$ by (\ref{id2.2}) it follows that $(V,a,b)$ is a subbimodule of $V.$ If $(V,a,b)=V,$ then $V=(V,a,b)=((V,a,b),a,b)=0$ by (\ref{id2.4}). Thus, $(V,a,b)=0.$

Consider $[V,a]$. Since $[\G,B]=0$  and $\cha B=3$, it follows from (\ref{id2.1}) that $[V,a]$ is a subbimodule of $V.$ Hence, $[V,a]=0$ or $[V,a]=V.$

Note that $D(a^3)=3a^2D(a)=0$. We will show next that if $D(b)=0$ for some $b\in \G$ then $[V,b]=0.$

Let $c,d\in \G.$ We have
\begin{eqnarray*}
(b,\bar{c},\bar{d})&=&\bar{bc}\cdot\bar{d}-b\cdot(\gamma cd + 2D(c)d + cD(d))\\
&=&\gamma bcd+2D(bc)d+bcD(d)-\gamma bcd-2bD(c)d-bcD(d)\\
	&=&2D(bc)d-2bD(c)d=0.
\end{eqnarray*}

It follows now from (\ref{id2.2}) that  $(V,b,x)$ is a subbimodule of $V$ for any $x\in \G$. If $(V,b,x)=V$ then $V=(V,b,x)=((V,b,x),b,x)=0,$ by superized (\ref{id2.5}). Thus $(V,b,\bar\G)=0.$

If $[V,b]=V$ then again from (\ref{id2.2}) we get
\bes
(V,B,B)=([V,b],B,B)\subseteq (BB,V,b)+B(B,V,b)+(B,V,b)B=0,
\ees
 which contradicts to the nonassociativity of $V.$ Hence, $[V,b]=0.$

Now, if $[V,a]=V,$ then $V=[V,a]=[[V,a],a]=[[[V,a],a],a]=[V,a^3]=0,$ by the above. Hence, $[V,a]=0.$
\ctd

\smallskip

The next proposition is an adaptation of the results from \cite{Sh2} to the bimodule situation.

\begin{Prop} \label{prop6.1}
 Let $V=V_{\0}\oplus V_{\1}$ be an irreducible bimodule over an alternative superalgebra $A=A_{\0}\oplus A_{\1}.$
If $N$ is a proper $A_{\0}$-subbimodule of $V_i$, and $(N,A_{\1},A_{\1})\subseteq N,$ then $V=(N\oplus A_{\1}NA_{\1})\oplus(NA_{\1}\oplus A_{\1}N).$
Moreover, if $[N,A_{\1}]_s=0$ then $N=V_i$.
\end{Prop}
\Proof It is clear that $(N+A_{\1}NA_{\1})\oplus(NA_{\1}+ A_{\1}N)$ is a proper subbimodule of $V$ and hence it is equal to $V$. Let us prove that the sums in the summands are direct.
Suppose that $N\cap A_{\1}NA_{\1}=L\neq 0.$ Then $L+L\circ A_{\1}$ is a nonzero subbimodule of $V.$ Hence, $L+L\circ A_{\1}=V,$ but this contradicts to the fact that $L$ is a proper $A_{\0}$-subbimodule of $V_i$. Therefore, $N\cap A_{\1}NA_{\1}=0.$

Now let $L=NA_{\1}\cap A_{\1}N.$ Then $LA_{\0}\subseteq NA_{\1}\cdot A_{\0}\subseteq NA_{\1}.$ On the other hand, $LA_{\0}\subseteq A_{\1}N\cdot A_{\0}\subseteq (A_{\1},N,A_{\0})+A_{\1}N\subseteq A_{\1}N.$ Hence, $LA_{\0}\subseteq L.$ Similarly, $A_{\0}L\subseteq L.$ Thus $L+L\circ A_{\1}=V$. But $LA_{\1}\subseteq NA_{\1}\cdot A_{\1}\subseteq N$ and, similarly, $A_{\1}L\subseteq N.$ Therefore $V=L+N,$ a contradiction. So,  $NA_{\1}\cap A_{\1}N=0.$

\smallskip
Finally, assume that $[N,A_{\1}]=0$. Then we have
\bes
(A_{\1}N)A_{\1}\subseteq (NA_{\1})A_{\1}\subseteq (N,A_{\1},A_{\1})+NA_{\0}\subseteq N,
\ees
therefore $V=N\oplus N\circ A_{\1}$ and $N=V_i$.
\ctd
\begin{Lem}\label{lem6.5}
 $(V_{\0},\G,\bar\G)=0$ or $(V_{\1},\G,\bar\G)=0.$
\end{Lem}
\Proof Let $W=\left\{v\in V\,|\,(v,\G,\bar\G)=0\right\}.$ We show first that $W\neq 0.$
Recall that $(\G,\G,\bar\G)=0$. By superization of (\ref{id2.4}), in view of Lemma \ref{lem6.4} we have
\bes
((V,\G,\bar\G),\G,\bar\G)&\subseteq &((\G,\G,\bar\G),V,\bar\G)+((V,\bar\G,\bar\G),\G,\G)\\&+&((\G,\bar\G,\bar\G),V,\G)+[V,\bar\G](\bar\G,\G,\G)\subseteq (V,\G,\G)=0.
\ees
Therefore, $W\neq 0$.   It follows from (\ref{id2.2})  that $W\G\subseteq W$.
Moreover, by (\ref{id2.5}) we have
\bes
((W,\bar\G,\bar\G),\G,\bar\G)\subseteq (W,\G,\bar\G)+(W,\bar\G,\G)=0.
\ees
Consequently, $(W,\bar\G,\bar\G)\subseteq W$. It is clear that $W$ is a graded subspace of $V$, and we see that the components $W_i,\ i=0,1,$ satisfy the condition of Proposition \ref{prop6.1}.
Consider $[W,\bar\G]$.  Since $\cha B=3$, it follows from (\ref{id2.1})  that $[W,\bar\G]$ is a $\G$-submodule of $V.$ By superization of (\ref{id2.6}) and (\ref{id2.7}),
\bes
([W,\bar\G],\bar\G,\bar\G)&\subseteq& ([\bar\G,\bar\G],\bar\G,W)+[(\bar\G,\bar\G,\bar\G),W]+ [(\bar\G,\bar\G,W),\bar\G]\\
&\subseteq& (\G,\bar\G,W)+[W,\bar\G]= [W,\bar\G].
\ees
Consequently, the subspaces $[W_i,\bar\G],\ i=0,1,$ satisfy the condition of Proposition \ref{prop6.1} as well.

Using (\ref{id2.1}) and Lemma \ref{lem6.4}, we obtain
$$
(\bar\G[W,\bar\G])\bar\G\subseteq ([W,\bar\G]\bar\G)\bar\G\subseteq [W,\bar\G]+[W,\bar\G]\G\subseteq [W,\bar\G].
$$
Therefore,  for every $i=0,1$, by Proposition \ref{prop6.1} we have
 $[W_i,\bar\G]=0$ or $V=[W_i,\bar\G] \bar\G+[W_i,\bar\G]$.

If $[W_i,\bar\G]=0$, then by Proposition \ref{prop6.1} $W_i=V_i$ and $(V_i,\G,\bar\G)=0$.

In the second case we have $V_{1-i}=[W_i,\bar\G]$. Then by (\ref{id2.8}), (\ref{id2.9}) we have
\bes
(V_{1-i},\bar\G,\G)&=&([W_i,\bar\G],\bar\G,\G)\subseteq [(\bar\G,\bar\G,\G),W]\\
&+& [(\bar\G,\bar\G,W),\G]+([\G,\bar\G],\bar\G,W)\subseteq [V,\G]=0,
\ees
hence $V_{1-i}\subseteq W$.
\ctd

\begin{Lem}\label{lem6.6}
If $(V_{\0},\G,\bar\G)=0,$ then $[V_{\0},\bar\G]=0.$
\end{Lem}
\Proof
Similarly as it was done for $W$  in the proof of Lemma \ref{lem6.5}, we can show that $[V_{\0},\bar\G]=0$ or $V=[V_{\0},\bar\G]\bar\G+[V_{\0},\bar\G].$
Consider the last case. As in the proof of Lemma \ref{lem6.5} we obtain that $(V_{\1},\G,\bar\G)=([V_{\0},\bar\G],\G,\bar\G)=0,$ so $(V,\G,\bar\G)=0.$

Note that for any nonzero $a\in \G$ we have $Va=V,$ hence $V=Va=(Va)a=Va^3.$ Thus $a^3\neq 0.$ Since  $\G$ is $D$-simple  \cite{Sh2}, it follows  from \cite{Yuan} that $\G$ is a field.

By (\ref{id2.2}) we have
$$
0=(V_{\0},\G,\bar\G)=([V_{\0},\bar\G]\bar\G,\G,\bar\G)=[V_{\0},\bar\G](\bar\G,\G,\bar\G).
$$
In particular, for any $a\in \G$ we have $0=[V_0,\bar\G](a,\bar{1},\bar{1})=2[V_0,\bar\G]D(a).$ (It was proved in \cite{Posner} that any $D$-simple algebra contains an identity element 1.) By choosing $a$ such that $D(a)\neq 0$ we obtain that $[V_0,\bar\G]=0.$
 \ctd

\smallskip

 Let $(A,D)$ be an algebra with a~derivation~$D$. An $A$-bimodule $(V,d)$ with a~linear mapping $d\colon V \to V$ is called \emph{a bimodule with~derivation}
or \emph{a $D$-bimodule} over the algebra $(A,D)$, if the linear
mapping $D +d\colon a+v \mapsto D(a)+d(v)$ is a derivation of the split null extension $E = A\oplus V$. If, in addition, $V$ has no proper
$d$-invariant $A$-subbimodules, then $(V, d)$ is called \emph{a $D$-simple} $A$-bimodule.

Let $(V, d)$ be an associative and commutative  $D$-bimodule over an algebra with~a derivation $(\G,D)$. Consider the split null extension $E = \G  \oplus V$
with the~derivation $D + d$ and construct the alternative superalgebra $B(E,D + d,\gamma)$. We have
$$
B(E,D + d,\gamma) = (\G \oplus  V ) \oplus \overline{\G \oplus V} =
(\G \oplus \bar\G) \oplus (V \oplus \bar V) = B \oplus W,
$$
where $W = V \oplus \bar V$ is an alternative bimodule over the superalgebra
$B= B(\G,D,\gamma)$. We denote this bimodule by $B(V, d,\gamma)$.

From the definition, it is easy to recover the explicit action of
$B(\G,D,\gamma)$ on $B(V, d,\gamma)$:
\begin{eqnarray*}
a \cdot v = v \cdot a = av, &
\bar a \cdot v = v \cdot \bar a = a \cdot \bar v = \bar v \cdot a =
\overline{av},
\\
\bar a \cdot \bar v = 2D(a)v + ad(v) + \gamma av, &
\bar v \cdot \bar a = 2d(v)a + vD(a) + \gamma va.
\end{eqnarray*}

We are now ready to prove

\begin{Th}\label{th6.2} Let $V$ be an irreducible faithful nonassociative  alternative bimodule over the simple superalgebra $B(\G,D,\gamma)$. Then, up to the changing of parity, $(V_{\0}, d)$ is a $D$-simple faithful
associative and commutative bimodule over the algebra
with derivation $(\G,D)$,
where $d(v) = -(v,\bar 1,\bar 1)$ is a derivation of $V_{\0};$ $V_{\1}=V_{\0}\cdot \bar{1}$ is an isomorphic copy of $V_{\0}$ and $V \cong B(V_{\0}, d,\gamma)$.
\end{Th}
\Proof By Lemmas \ref{lem6.4}, \ref{lem6.5}, \ref{lem6.6}, we may assume, up to the changing of the parity, that $(V_{\0},\G,B)=[V_{\0},B]=0$.  Let $d:V_{\0}\rightarrow V_{\0},\ d: v \mapsto 2(v,\bar 1,\bar 1)$, then it follows from (\ref{id2.2}) that $(V_{\0},d)$ is an associative and commutative $D$-bimodule over $(\G,D)$.
For $v\in V_{\0}$, denote $\bar v=v\cdot \bar 1$, then we have
\bes
a\bar v&=&a(\bar 1 v)=(a\bar 1)v=\bar a v=(\bar 1 a)v=\bar 1(av)=\overline{av},\\
\bar a\bar v&=&\bar a(\bar 1 v)=(\bar a\bar 1)v-(\bar a,\bar 1,v)=(2D(a)+\g a)v-(a\bar 1,\bar 1,v)\\
&\stackrel{(\ref{id2.6})}=&(2D(a)+\g a)v-(\bar 1\cdot\bar 1,a,v)-(\bar 1,\bar 1,v)a-
(\bar 1,a,v)\bar 1\\
&=&2D(a)v+d(v)a+\g av,
\ees
and similarly $\bar v\bar a=2d(v)a+vD(a)+\g	va$. Therefore, $V\cong B(V_{\0},d,\g)$.
\ctd

\smallskip
The converse is also true:

\begin{Prop}\label{prop6.2}
For any $D$-simple faithful associative and commutative
bimodule $(V, d)$ over $(\G,D)$ and any $\gamma \in \G$
 the bimodule $B(V, d,\gamma)$ is  an irreducible faithful alternative superbimodule over the simple superalgebra $B=B(\G,D,\gamma)$. Moreover, the $B$-superbimodules $B(V, d,\gamma)$ and $B(V', d',\gamma)$
are isomorphic if and only if the differential bimodules $(V, d)$ and $(V', d')$ over $(\G,D)$ are isomorphic. 
\end{Prop}

\Proof Let $W$ be a nonzero subbimodule of $V$, then $W=W_{\0}+\bar U$ for some $\G$-subbimodule $U$ of $V_{\0}$. We show that $W_{\0}\neq 0.$ Assume on the contrary that $W_{\0}=0.$ Then  $\bar{v}\cdot\bar{1}=\bar{1}\cdot\bar{v}=0$ for any $\bar{v}\in W_{\1}=\overline{U}.$ Hence $d(v)=0$ for any $v\in U.$ It follows that $U$ is a nonzero $D$-submodule of $V_{\0}.$ Thus $U=V_{\0}$ and $0=d(va)=vD(a)$ for any $v\in V_{\0},$ $a\in\G.$ Consequently, $D(\G)=0,$ a contradiction.

Therefore $W_{\0}$ is a nonzero $\G$-subbimodule of $V_{\0}$ and $\overline{W_{\0}}\subseteq W_{\1}.$  For any $v\in W_{\0}$, $\bar{1}\cdot\bar{v}=\g v+d(v)\in W_{\0}.$ Hence $d(W_{\0})\subseteq W_{\0}$ and $W_{\0}=V_{\0}.$

Finally, if the superbimodules $B(V, d,\gamma)$ and $B(V', d',\gamma)$ are isomorphic then their even parts $V_{\0}$ and $V'_{\0}$ are isomorphic associative and commutative bimodules over $\G$. Moreover, the derivations $d$ and $d'$ are defined via the same mapping $v\mapsto 2(v,\1,\1)$, hence $(V,d)\cong (V',d')$.
\ctd

\medskip
Clearly, for any $D$-simple algebra $(\G,D)$ the regular bimodule with derivation $\Reg (\G,D)$ is a $D$-simple bimodule over $(\G,D)$. Moreover, for any $a\in \G$  the map $D+R_a: \G\rightarrow \G$ is a derivation of the bimodule $\Reg\G$. In fact, for any $v\in\Reg\G,\ b\in\G$ we have
\bes
(D+R_a)(vb)=D(vb)+vba=D(v)b+vD(b)+vba=(D+R_a)(v)b+vD(b).
\ees
It is clear that the differential bimodule $(\Reg \G,D+R_a)$ is $D$-simple as well.

\begin{Prop}\label{prop6.3}
The $D$-bimodules $(\Reg \G,D+R_a)$ and $(\Reg \G,D+R_b)$ over $(\G,D)$ are isomorphic if and only if there exists an invertible $c\in\G$ such that $D(c)=(a-b)c$. In particular, if $a,b\in F$ then $(\Reg \G,D+R_a)\cong( \Reg \G,D+R_b)$ if and only if $a-b\in Spec\,D$.
\end{Prop}
\Proof Let $\f:(\Reg \G,D+R_a)\rightarrow (\Reg \G,D+R_b)$ be an isomorphism of $D$-bimodules over $(\G,D)$. Let $c=\f(1)$, then for any $x\in\G$ we have $\f(x)=\f(1\cdot x)=\f(1)x=cx$. Since $\f$ is an isomorphism, $c\G=\G$ and $c$ is invertible. Futhermore, 
\bes
\f((D+R_a)(x))&=&cD(x)+cax=(D+R_b)(\f(x))=D(cx)+bcx\\
&=&D(c)x+cD(x)+bcx,
\ees   
which gives $(D(c)-(a-b)c)\G=0$. Since $\G$ is $D$-simple, this proves that $D(c)=(a-b)c$. Conversly, if $c$ is an invertible element that satisfies this relation then the application $x\mapsto cx$ is an isomorphism of $D$-bimodules $(\Reg \G,D+R_a)$ and $(\Reg \G,D+R_b)$.
 It remains to notice that if $\l\in Spec\, D$, that is, $D(c)=\l c$ for some $c\in\G$ then $c\G$ is a $D$-ideal in $\G$, hence $c\G=\G$ and $c$ is invertible.
\ctd

\smallskip

We now describe the faithful  $D$-simple bimodules over a finite dimensional $D$-simple associative commutative algebra $(\G,D)$ in the case when the ground field $F$ is algebraically closed. Observe that due to \cite{Yuan}   in this case $\G\cong F[t_1, \ldots,t_n]/(t_1^3,\ldots ,t_n^3)$.

\begin{Th}\label{th6.3}
Let $(\G,D)$ be  a  $D$-simple associative commutative algebra over an algebraically closed field $F$ of characteristic $3$, and let  $(V,d)$ be a  $D$-simple finite dimensional associative and commutative bimodule  over $(\G,D)$. Then there exist $\l\in F$ and $v\in V$ such that $V=v\G$ and the action of the derivation $d$  is given by $d(va)=v(\l a+D(a))$. In other words, the bimodule $(v\G,d)$ is isomorphic to the  $D$-bimodule $(\Reg \G,D+\l)$. The parameter $\l$ is defined uniquely modulo $Spec\,D$.

\end{Th}
\Proof 
Since the field $F$ is algebraically closed and $V$ is finite dimensional, the linear operator $d:V\rightarrow V$ has an eigen  vector $v$ with an eigen value $\l\in F$. For any $a\in\G$ we have $d(va)=d(v)a+vD(a)=v(\l a+D(a))$, which proves that $v\G$ is a $d$-subbimodule of $V$. Clearly,
$v\G\neq 0$ and therefore $V=v\G$.

Consider the mapping $\phi : \Reg\G\rightarrow V,\ a\mapsto va$. Clearly, $\phi$ is an isomorphism of $\G$-modules. Moreovere, we have
\bes
\phi((D+\l)(a))=v(D(a)+\l a)=d(va)=d(\phi(a)),
\ees
which shows that $(V,d)\cong (\Reg\G,D+\l)$. The last statement follows from Proposition \ref{prop6.3}.
\ctd

As a~result, for the finite dimensional alternative superalgebras over 
an algebraically closed field we have a complete description of irreducible 
nontrivial bimodules, which are not associative.

\begin{Th}\label{th6.4}
Let $A$ be a finite dimensional alternative superalgebra 
over an algebraically closed field~$F$, $V$ be an irreducible 
faithful nontrivial alternative superbimodule over~$A,$ 
which is not associative. Then either $\dim A\leq 2,\ A_{\1}=Fx,\, x^2=0$ 
and $V\cong V^{\e}(1|1),\, \e^2+\e+1=0$,
or   the superalgebra~$A$ is simple and $V$ is either isomorphic to one of the bimodules $\Reg A$, $(\Reg A)^{\mathrm{op}}$,
or one of the following cases holds:
\begin{itemize}
\item
$A = F[\sqrt{1}\,]$ and $V$ is one of the two-dimensional bimodules 
of types \textup{1) -- 3)} described in~Theorem~\ref{th2.2};
\item
$\cha F = 3$, $A = B(1|2)$, $V = V_{\lambda,\mu}(3|3)$;
\item
$\cha F = 3$, $A = B(\G,D,\g)$, $V = B(\Reg\G, D+\l,\g),\l\not\in Spec\,D$;
\item
$\cha F = 2$, $A = \mathrm M_{1|1}(F),\ V=\Cay L$,  $L$ is a minimal $\Z_2$-graded left ideal of $\mathrm M_{1|1}(F)$;
\item $\cha F = 2$, $A = \mathrm M_2(F)[\sqrt1],\ V=(\Cay L)[\sqrt1]$,  $L$ is a minimal left ideal of $\mathrm M_2(F).$
\end{itemize}
\end{Th}

Observe that there is a misprint in the announcement\textit{} of this result in \cite{TSh}:
the case $A=B(\G,D,\g)$ is omitted.

\end{document}